\documentclass[a4paper,10pt]{amsart}
\usepackage{mathrsfs}
\usepackage{calrsfs}
\usepackage{times}
\usepackage[scaled=0.92]{helvet}
\usepackage[greek,english]{babel}
\usepackage[iso-8859-7]{inputenc}

\newcommand{\lc}{\left\lceil}
\newcommand{\rc}{\right\rceil}
\newcommand{\lf}{\left\lfloor}
\newcommand{\rf}{\right\rfloor}

\usepackage{amsfonts}
\usepackage{amssymb}
\usepackage[usenames]{color}
\usepackage[active]{srcltx}
\newtheorem{theorem}{Theorem}   
\newtheorem{conj}[theorem]{Conjecture}
\newtheorem{lemma}[theorem]{Lemma}
\newtheorem{proposition}[theorem]{Proposition}

\newtheorem{remark}[theorem]{Remark}

\title{On holomorphic polydifferentials in positive characteristic }
\author{Sotiris Karanikolopoulos}
\address{
Department of Mathematics, University of the \AE gean, 83200 Karlovassi, Samos,
Greece
}
\email{mathm03005@aegean.gr}

\date{\today}

\begin{document}
\bibliographystyle{amsplain}

\begin{abstract}
In this paper we study  the space $\Omega(m)$,
of holomorphic $m$-(poly)differentials of a function field of a curve defined over an algebraically closed field of characteristic $p>0$ when
 $G$ is cyclic or elementary abelian group of order $p^n$; we give bases for each case when the base field is rational, introduce the 
Boseck invariants and give an elementary approach to the $G$ module structure of $\Omega(m)$  in terms of Boseck invariants. The last 
computation is achieved without any restriction on the base field in the cyclic case, while in the elementary abelian case it is assumed that 
the base field is rational. An application to the computation of the tangent space of the deformation 
functor of curves with automorphisms is given.
\end{abstract} %\footnote{abstract and intro need fill in}

 \thanks{{\bf keywords:} Automorphisms, Curves, Differentials, Galois module
 structure, {\bf AMS subject classification} 14H37,11G20}
\maketitle

\section{Introduction}
Let $F$ be an algebraic function field with field of constants  $K$, where $K$ is 
 an algebraic closed field of characteristic $p$.
Let $F/E$ be
a Galois extension with abelian Galois group $G$ of order $p^n$.
We will denote by $\Omega_{F}(m)$ the space of  holomorphic $m$-(poly)differentials of $F$.
We know that $\Omega_{F}:=\Omega_F(1)$ is a $g_F$-dimensional $K$-space, while
the $\Omega_{F}(m)$, is a $(2m-1)(g_F -1)$-dimensional $K$-space, when $g_F\geq
2$.
%  An $m$-holomorphic differential $\omega$ is an assignment of a holomorphic
%function to each local coordinate $z$ on $F$ so that $f(z)(dz)^{\otimes m}$ is
%invariantly defined. We see them as the $mth$ powers of holomorphic ($m=1$)
%differentials. 
The Galois module structure of the space of holomorphic $1$-differentials has been determined explicitly in some cases.
The cyclic group case was studied by  Hurwitz \cite{Hu:1893} if the characteristic of $K$ is zero.
When  $F/E$ is unramified and $G$ has a prime to $p$ order, or is a cyclic
group, Tamagawa \cite{Tamagawa:51} proved that is the direct sum of one identity
representation of degree one and $g_E -1$ regular representations. Valentini
\cite{Val:82} generalized this result for unramified extensions having
$p$-groups as their Galois groups, while Salvador and Bautista
\cite{Salvador:00} determined completely the semisimple part of holomorphic
differentials when  $G$ is a $p$-group. If $G$ is cyclic then Valentini and
Madan \cite{vm}, determine completely the whole structure of $\Omega_{F}$ in
terms of indecomposable $K[G]$-modules. The same is done when $G$ is an
elementary abelian, by Calder{\'o}n, Salvador and Madan \cite{csm}.
Also N. Borne \cite{borne06} developed a theory, using advanced techniques from 
modular representation theory and $K$-theory, and he is able to compute the 
$K[G]$-module structure for holomorphic 
$m$-differentials $\Omega_F(m)$, when $G$ is a cyclic group of order $p^n$.
In general the $K[G]$-module structure of $\Omega_{F}$, in positive
characteristic is unknown. The difficulties that arise in positive 
characteristic, in contrast to the same problem in characteristic zero,
are first all the difficulties of modular representation theory, in contrast to 
ordinary representation theory, and second the appearance of wild ramification 
in extensions $F/E$.

In this paper we will mainly focus on the $m>1$ case and on the two ``extreme'' cases of 
abelian groups of order $p^n$ namely cyclic groups and elementary abelian groups.
We first compute a basis for holomorphic differentials and then we define 
the quantities  $ \nu_{ik}(m)$ to be
$$
\left\lfloor\frac{m\delta_i + \{ \textrm{evaluation of the kth  E- basis element of F  by a normalized valuation of F }\} }{p^{e_i}}\right\rfloor,
$$
where $\delta_i, e_i, i$ are related with the ramification of the extension, see Conjecture \ref{conj} below and Remark \ref{unifying remark}; the basis element is evaluated by a (normalized) valuation determined by a place of $F$ above a ramified place of $E$ and $\lfloor \cdot \rfloor$ denotes the integer part. 
We introduce the {\em Boseck invariants} 
that are quantities of the form $\Gamma_k(m):=\sum_{i} \nu_{ik}(m)$, where 
the index $i$ runs over all ramified places.
These quantities were used by Boseck \cite{Boseck}  for constructing bases for
1-holomorphic differentials and the conditions for holomorphicity 
are expressed in terms of them. In order to find Bosek's invariants for the function fields $F/E$ that we study, we take the rational
 extensions $F/K(x)$ and find $K$-bases for the corresponding $\Omega_F(m)$'s for $m\geq1$ 
(when $m=1$ these bases often called Boseck bases). The choice of the rational function field is clear each time from the defining 
equations of our curves (Eq. (\ref{eq.cyclic}),  (\ref{kumeq}), (\ref{equation e.abelian})).
It turns out that in the cases we study, the Boseck invariants determine the
Galois module structure
of the space of $m$-holomorphic differentials. 
{  The formula that gives this $K[G]$ module structure in terms of the Boseck invariants remains the same in both the elementary
abelian and the cyclic case. In addition in the cyclic case this formula is independent from $r$, 
with $r$ measuring  how high in the 
tower of intermediate fields is placed an unramified subextension $E_r /E$. This  is not true when $m=1$ (see Remark \ref{m=1}).}
% This structure seems to be unaffected by the ramification type ($F/E$ is totally
% ramified, that is the elementary abelian case, or not, the cyclic case; in the last case we also observe an i) when $m\gneqq1$. 
% while for the $m=1$ case the dependence from
% the ramification type is inevitable.
%  However, one can use Tamagawa and Valentini
% results to treat the unramified $m=1$ case (or the unramified subextension that
% appears  when no place is totally ramified in the cyclic case, as Valentini and
% Madan did in \cite{vm}). 
The elements $\Gamma_{k}(m)$ carry a lot of information, 
for instance the  degree of the different of the extension each time,
 can be described  totally by them (see Remarks \ref{degdif1}, \ref{degdifKummer},
\ref{degdif3}), namely
$$\frac{2}{2m-1}\sum_{k}\Gamma_k(m) =\deg\mathrm{Diff}(F/E).$$
%  This problem is connected to the study of the tangent space of the global deformation functor of curves with automorphisms (see \cite{KoJPAA06} and  \cite{KoArxiv2}).  
% Some of the results for the $m=1$ case have remarkable similarities. That was the reason for looking  some quantities who would help to describe all the $K[G]$-module structures of $\Omega_F(m)$ for $G$ with different structures and $m$'s. Here we focus on the $K[G]$-module structure of $\Omega_F(m)$, for $m>1$ and see that case as a generalization of the $m=1$ case. We introduce the \textit{Boseck invariants}. These are quantitities that describe completely the $K[G]$-module structure of holomorphic 
% $m$-differentials , for $m\geq1$. 
% 
% . Boseck invariants are quite powerfull in the sense that,  the degree of the different of the extension each time, can be described  totally by them (see Remarks \ref{degdif1}, \ref{degdif2}, \ref{degdif3}). That explains their invariance of the base field ( rational or not ). They also describe the $K[G]$-module structure of holomorphic of $\Omega_F(m)$, for every $m\geq 1$.

Using Boseck invariants, we describe every $K[G]$-module structure of holomorphic 
$m$-(poly)differentials $\Omega_F(m)$ for $m\gneq1$, when: $G$ is cyclic (for $m=1$ case this is \cite[Theorem 2]{vm}), or elementary abelian of degree $p^n$ (for $m=1$ case this is \cite[Theorem 1]{csm}). 
Finally we show how the Boseck basis in the $m=1$ case   when  $F/E$ is tame, i.e,  of degree $n$, with $n$ being prime to $p$, will give another proof of the classical 
 Hurwitz Theorem \cite[p. 439, formula 2]{Hu:1893}.

Our approach, which is quite elementary, follows closely the ideas of Valentini-
Madan,
 Calder{\'o}n,Villa-Salvador  and Madden \cite{Madden78}. Madden, used the same
analysis and constructed a $K$-basis of  $\Omega_F (1)$,  in order to compute
the rank of the Hasse-Witt matrix. We should mention here that all the above
authors had used Boseck invariants in their papers \cite{vm}, \cite{csm},
\cite{Madden78}.

%For simplicity, we focus in two extreme, for the structure of $G$, cases : namely when $G$ is cyclic, of order $p^n$ and in contrary, when $G$ is elementary abelian, of order $p^n$, with $E=K(x)$.

The organization of the paper is as follows: In first  section we focus on the
cyclic case and give the Galois module structure of $\Omega_F(m)$ in Theorem \ref{basic thrm}, in subsection \ref{Galoismodulecyclic}.
As our analysis is going deeper, the ideas of Boseck \cite{Boseck} are rising
up. We follow him, in subsection \ref{cyclic basis}, in order to give a basis
for $m$-holomorphic (poly)differentials of $F$, when $E$ is rational.
Subsection  \ref{HuClass} is devoted to a proof of the classical result of Hurwitz mentioned 
above.
In section \ref{elab}  we consider the elementary abelian
case. At first we give an analogous $K$ basis for $\Omega_{F}(m)$. These bases
can lead to the computation of $m$-Weierstrass points (see Remark \ref{wp}).
%  and
% they also are of crucial importance, because Boseck invariants arise from them. 
Our proof of theorem \ref{fff} in  subsection \ref{elabgaloisstructure} is based on the work 
 of Calder{\'o}n, Salvador and Madan \cite{csm}.

In section \ref{consection} we state a conjecture concerning the Galois module 
structure  of $m$- differentials when $G$ is a  general abelian extension 
of order $p^n$.

Finally, in section \ref{4section} we give an application to
the computation of the tangent space of the local deformation functor in the
sense of J.Bertin and A. M\'ezard \cite{Be-Me}. The results are given in terms 
of the Boseck invariants, and coincide with computations done previously by
other authors \cite{Be-Me},\cite{KontoANT},\cite{KoJPAA06},\cite{KoArxiv2}
using completely different methods. 
This allows us to verify our complicated computations concerning 
the Galois module structure.
Also proving the conjectures stated in section \ref{consection} will give  a
method in order to 
compute the above mentioned dimension in the case of abelian groups, a problem
that is still open. 

An other application we have in mind and we would like to explore in a following
article is  the computation of higher order Weierstrass points and the study of
the fields generated by the coefficients of them, a problem that is similar to
the classical study of fields generated by torsion points of the Jacobian,
\cite{Silverman92}, 
\cite{Towse96}. 

In order to avoid trivial cases we will always assume that $g_F \geq2$, where
$g_F$
 is the genus of $F$. We use the symbol $\mathbb{P}$ to denote the set of places
(sometimes referred just as ``primes'') of the field in question.

\section{The Cyclic Case}

Since  the characteristic $p$ divides the order of the group $G$,    the representation of $G$ on 
$\Omega_F(m)$
 is not necessary completely reducible, but it is the direct sum of indecomposable $K[G]$-modules.
Let $\sigma$ be a generator for $G$. 
The unique indecomposable $K[G]$-module of degree $k$ is isomorphic to 
$K[G]/\left\langle (\sigma -1)^k \right\rangle $ \cite[p.156, Ex. 1.1 ]{Wein}.
 For $k=1$ we obtain  the identity representation
and for $k=p^n$ the regular representation. 
% Let $d_k$ be the number of times that the above indecomposable
% representation of degree $k$, occurs in  $\Omega_{F}(m)$ seen as a $K[G]$-module.

Let 
\begin{equation} \label{dec}
 \Omega_{F}(m):= \bigoplus_{\lambda =1}^t  W_\lambda,
\end{equation}
 be a decomposition into a direct sum
of indecomposable $K[G]$-modules and let $d_k$ be the number of $W_{\lambda}$'s that are isomorphic to
$K[G]/\left\langle (\sigma -1)^k \right\rangle $. We will  compute the $d_k$'s. First we define
\[
 \Omega_{F}^{i}(m)=\{\omega \in \Omega_{F}(m) : (\sigma -1)^i \omega =0 \}, \mbox{ for } i=0, 1,\ldots p^n .
\]
These $K$-subspaces form an increasing sequence with $\Omega_{F}^{0}(m)=0$ and
$\Omega_{F}^{p^n}(m)=\Omega_{F}(m)$, while
\begin{equation} \label{dec2}
 \dim_K \Omega_{F}^{i}(m)=\sum_{\lambda =1}^t \dim_K (W_\lambda \cap \Omega_{F}^{i}(m)).
\end{equation}

%The quotients of those spaces play a major role in the calculation of $d_k$'s.

\begin{lemma} \label{lem1}
We have $\dim_K K[G]/\langle (\sigma-1)^k \rangle=k$.
\end{lemma}
\begin{proof}
% Observe that  $\dim_K \mathrm{ker}(\sigma-1)^k=k$ 
%therefore $\dim_K \mathrm{im}(\sigma-1)^k=p^n-k$ and the quotient 
%$\mathrm{dim}_k K[G]/\langle (\sigma-1)^k \rangle=\mathrm{dim}_k K[G]/\langle \mathrm{im}(\sigma-1)^k \rangle$ 
%has dimension $p^n-(p^n-k)=k$.
%{\bf 8 Sep 08 Write more here}
Consider the map $(\sigma-1)^k:K[G] \rightarrow K[G]$.
Observe that   $\dim_K \mathrm{ker}(\sigma-1)^k=k$ therefore $\dim_K \mathrm{im}(\sigma-1)^k=p^n-k$ and the quotient 
$\mathrm{dim}_K K[G]/\langle (\sigma-1)^k \rangle=\mathrm{dim}_K K[G]/\langle \mathrm{im}(\sigma-1)^k \rangle$ 
has dimension $p^n-(p^n-k)=k$.

% Let $\phi:K[G] \rightarrow K[G]/ \left\langle (\sigma -1)^k \right\rangle$ be the canonical homomorphism of vector spaces over $K$. Then 
% \begin{eqnarray}
% \dim_K (K[G]) &=& \dim_K (\ker \phi) +\dim_K (\textrm{im} \phi)\nonumber \\
% \dim_K (K[G]) &=& \dim_K (\left\langle (\sigma -1)^k \right\rangle)+\dim_K K[G]/\langle (\sigma-1)^k \rangle \nonumber\\
% p^n &=& p^n-k + \dim_K K[G]/\langle (\sigma-1)^k \rangle  \nonumber                                                     \end{eqnarray}
\end{proof}

\begin{lemma} \label{lem2}
$$\dim _K \{\alpha\in K[G]/\left\langle (\sigma -1)^k \right\rangle : (\sigma -1)^i \alpha=0\}=\begin{cases}  i, & \textrm{if }  i\leq k ,\cr
              k, & \textrm{if } i \gneqq k. \cr  \end{cases}$$
\end{lemma}
\begin{proof}
Indeed, we would like to compute  the kernel of the multiplication with
$(\sigma -1)^i $, $$(\sigma -1)^i :K[G]/\left\langle (\sigma -1)^k \right\rangle \rightarrow K[G]/\left\langle (\sigma -1)^k \right\rangle.$$
We distinguish the following two cases:
\begin{enumerate}
\item
If $i\leq k$ then 
$\ker(\sigma -1)^i = (\sigma -1)^{k-i} K[G]/\left\langle (\sigma -1)^k \right\rangle$ and
$\dim_K (\sigma -1)^{k-i} K[G]/\left\langle (\sigma -1)^k \right\rangle=p^n -(k-i)-(p^n -k)=i$.

 \item 
If  $i\geq k$ then  $\ker(\sigma -1)^i =K[G]/\left\langle (\sigma -1)^k \right\rangle$ and
$\dim_K  K[G]/\left\langle (\sigma -1)^k \right\rangle =k$ according to Lemma \ref{lem1}.

\end{enumerate}
\end{proof}
\begin{sloppypar}
Using the  decomposition of  $\Omega_{F}(m)$ given in Eq. (\ref{dec2}) and Lemma \ref{lem2} we obtain
$\dim_K \Omega_{F}^{i}(m) =\sum_{k=1}^{i-1} kd_k  +\sum_{k=i}^{p^n} id_k .$ So 
\end{sloppypar}
\begin{equation}
\dim_K  \left( \Omega_{F}^{i+1}(m)/\Omega_{F}^{i}(m)\right)=\sum_{k=i+1}^{p^n }d_k.
\end{equation}
Therefore, for $k=1,\ldots,p^n -1$ we have:
\begin{eqnarray}\label{dks}
d_{p^n} &=& \dim_K  \left(\Omega_{F}^{p^{n}}(m)/\Omega_{F}^{p^{n}-1}(m)\right), \\
d_k &=& \dim_K  \left(\Omega_{F}^{k}(m)/\Omega_{F}^{k-1}(m)\right)-\dim_K  \left(\Omega_{F}^{k+1}(m)/\Omega_{F}^{k}(m)\right).\nonumber
\end{eqnarray}
Following \cite{vm} we write down a convenient $E$-basis for $F$ and we find the $G$-action on the basis elements and state the main Theorem. 

Since $F/E$ is cyclic of degree $p^n$, there is a tower of intermediate fields
\begin{equation} 
\label{Artin-inter} 
E=E_0 \subset E_1 \subset E_2 \subset \cdots \subset E_n =F,
\end{equation}
where each of the $E_{j}/E_{j-1} $
is an  Artin-Schreier extension given by 
\begin{equation}\label{eq.cyclic}
E_{j}=E_{j-1}(y_j ),\; y_{j}^p -y_j =b_j ,\; b_j \in E_{j-1},\; 1\leq j\leq n.
\end{equation}
The elements   $b_j $  are called to be in standard form, for a place $P$ of $E_{j-1}$, if the valuation of the divisor of $b_j$ at $P$  is positive, zero or relatively prime to the characteristic  $p$.

The first part of  the following Theorem is due to Madden, and allows us to take
 $b_j $'s in standard form \cite[Theorem 2, p. 308]{Madden78}, while the other
part is due to R. Valentini and M. Madan \cite[Lemma 1, p. 109]{vm} and allows us  to
select a convenient 
$E$-basis of $F$.

\begin{theorem}\label{E-basis}
 The elements $y_j$ and $b_j$ can be selected so that:
\begin{enumerate}
 \item For any place $P$ of $E_{j-1}$ divisible by a ramified place in $F/E$ the valuation of $P$ of the 
divisor of $b_j$ is either zero or negative and relatively prime to $p$.
\item $\sigma^{p^{j-1}}(y_j)=y_j+1$.
\end{enumerate}
For $0 \leq k \leq p^n-1$ consider its $p$-adic expansion
$
 k:=a_1^{(k)}+a_2^{(k)} p+ \cdots + a_n^{(k)} p^{n-1}
$
and denote by $w_k=y_1^{a_1^{(k)}}y_2^{a_2^{(k)}}\cdots y_n^{a_n^{(k)}}$.
Then $F$ is an $E$ vector space with basis $ \{w_k : 0\leq k\leq p^n -1 \}  $. The
$G$-action on the $w_k $'s is given by
$$(\sigma -1)^k w_k=\prod_{\epsilon=1}^n a_{\epsilon}^{(k)}!$$
\end{theorem}
This basis, has the following property
\begin{lemma} \label{old2}
Let $\bar{P}$ be a place of $E$ and let $P_1 , P_2 ,\ldots, P_r $ be the places of $F$,
above $\bar{P}$. Let $v_i$ the normalized valuation of $F$, determined by
$P_i$, $i=1,\ldots, r$. Let also $b_j$ be in standard form for any place of
$E_{j-1}$ below some $P_i$. If $z=\sum_{k=0}^{p^{n}-1}c_k w_k $, then $\min_i v_i (z)=\min_{i,k}v_i (c_k  w_k )$.
\end{lemma} 
\begin{proof}  \cite[Lemma 2, p.109]{vm} or \cite[Lemma  3, p.310]{Madden78} .
\end{proof}

\subsection{Galois Module Structure of $\mathbf{ \Omega_F (m)}$}\label{Galoismodulecyclic}

Let  $\bar{P_{i}}, i=1,\ldots ,s $  be the places of $E$ which ramify in $F$ and 
set  $p^{e_{i}}:=e(P_F /\bar{P_{i}})$, $i=1,\ldots ,s$ for the corresponding ramification 
indices, with $P_F $ a fixed place of $F$ above $\bar{P_{i}}$.  We will denote by  $\bar{v}_i $  the
normalized valuation of $E$ determined by $\bar{P}_i $. Set $r=n -\max_{i}\{e_i \}$. We observe that $E_r / E$ is an unramified extension:
% We can consider the decomposition field $T\subset E_r $. Then $F/E_r $ is a
% totally ramified extension and in our case $E_r /E$ will be unramified because of our set up:
if not then from the transitivity of the ramification indices we will have
$$e(P_F /\bar{P_i } )=e(P_F /P_{E_{r}})\cdot e(P_{E_{r}}/\bar{P_i } )
\Rightarrow p^{e_i }=[F:E_r ]\cdot p^h \Rightarrow p^{e_i }=p^{n-r}\cdot p^h ,$$
where $h$ is a nonzero natural number. So $e_i =n-r+h \stackrel{r=n-\max_{i}\{e_i \}}{\Longrightarrow} e_i =\max_{i}\{{e_i }\}+h$ which
contradicts the maximality of $e_i .$ 

Fix an $i$. Let $r_i =n-e_i$. Let also $P(i,j,\mu)$ be the places of $E_j$
which divide $\bar{P}_i $ and $v(i,j,\mu,z)$ be the normalized valuation of  $ E_j $ determined
by $P(i,j,\mu)$, applied to an element $z\in E_j$. Each of the $E_j /E_{j-1}$ is normal and
separable, so every one of the $P(i,j,\mu)$'s will have the same exponent
$d(P(i,n,\mu)/\bar{P}_i ):=\delta_i$ in the different of $F/E$, $\mathrm{Diff}(F/E)$.
We can recover this different from the $\mathrm{Diff}(E_j /E_{j-1}),$ for all $\; j=1,\ldots,n$ using
the transitivity property of the different (see Stichtenoth, \cite[p.88, Corollary III.4.11.(a)]{StiBo}) $$\mathrm{Diff}(E_j /E)=\mathrm{Con}_{E_j /E_{j-1}}(\mathrm{Diff}(E_{j-1} /E))+\mathrm{Diff}(E_j /E_{j-1}).$$
Also, every automorphism of $F$
will act transitively on every place over $\bar{P}_i$, so the set $\mathrm{Diff}(E_j /E_{j-1})$
is stable under $\sigma$, and the exponent of $P(i,j,\mu)$ in $\mathrm{Diff}(E_j /E_{j-1})$ is independent of $\mu$.
The following Lemma gives us the relation among $\delta_i$ and $b_j$.
\begin{lemma}\label{different} %\footnote{i gave a sketch of proof}
 Let $\Phi(i,j)=-v(i,j-1,\mu,b_j)$. The different $\mathrm{Diff}(F/E)$ is given by
\[
 \mathrm{Diff}(F/E)=\sum_{i=1}^s \delta_i \sum_\mu P(i,n,\mu),
\]
where 
\[
 \delta_i=(p-1) \sum_{j=n-e_i+1}^n (\Phi(i,j)+1) p^{n-j},
\]
which  equals to
\[
 \delta_i=(p-1) \sum_{j=1}^n \Phi(i,j) p^{n-j}+(p^{e_i}-1).
\]

The valuations of the basis elements $w_k$ are given by 
\[
v(i,n,\mu,w_k)=-\sum_{j=1}^n a_j^{(k)} \Phi(i,j) p^{n-j}. 
\]
\end{lemma}
\begin{proof}This comes from Proposition 2 of Madden \cite{Madden78}, after replacing $j$, by $\nu+n-e_i$ and noticing that the different exponent does not depend on the choice of the base field (rational or not). However, we can prove it directly; fix an $i$ and apply the transitivity property of the different exponent  $d(P(i,n,\mu)/\bar{P}_i )$, which equals to:
$$e\left(P(i,n,\mu)/P(i,n-1,\mu)\right) d\left(P(i,n-1,\mu)/\bar{P}_i \right)+d\left(P(i,n,\mu)/P(i,n-1,\mu)\right),$$  $n-r_i=e_i$ times to get $\delta_i$. Observe finally that  $v(i,j-1,\mu,b_j)=0$, for all $1\leq j \leq n-e_i$. The last equality comes from the  $E$-basis of Theorem \ref{E-basis} and the fact that: $-\Phi(i,j)=v(i,j-1,\mu,b_j)=v(i,j,\mu ,y_j)$. Now, using the transitivity of the valuations, one can easily get that $v(i,n,\mu,y_j)=-\Phi(i,j)p^{n-j}$. 

\end{proof}
We are now  ready to define the key-quantities for our Theorem. For $k=0,1,\ldots,p^n -1,$ we define
\begin{equation}\label{niks}
\nu_{ik}(m):=\biggl\lfloor \frac{m\delta_i - \sum_{j=1}^{n}a_j ^{(k)}\Phi(i,j)p^{n-j}}{p^{e_i }}\biggr\rfloor,
\end{equation}
and $\Gamma_k(m):=\sum_{i=1}^s \nu_{ik}(m)$. In subsection \ref{cyclic basis}  we  interpret these $\Gamma_k(m)$ 's  as the Boseck invariants. %\footnote{i added this sentence}

%\footnote{i added this for reader's convenience} 
\begin{theorem}\label{basic thrm}

Let $G$ be a cyclic group of automorphisms of $F$, with $|G|=p^n .$ 
Set  $E=F^G$ and let  $g_E$ be the genus of $E$. Let  $m$ be  a natural number with $m>1$. 
The regular representation of $G$ occurs $d_{p^n}=\Gamma_{p^{n}-1}(m)+(g_E -1)(2m-1)$   times in the representation of $G$ on
 $\Omega_F (m)$. For $k=1,\ldots, p^n -1$, the indecomposable representation of degree $k$ occurs 
$d_{k}=\Gamma_{k-1}(m) -\Gamma_{k}(m)$ times.\\

\end{theorem}
\begin{proof}

As we saw, we have to compute $\dim_K (\Omega_{F}^{k+1}(m)/\Omega_{F}^{k}(m)),$ for $k=1,\ldots, p^n -1.$ 
Choose an $x \in E$,  such  that $dx \neq 0$. Every holomorphic (poly)differential $\omega $ of 
$F$ can be written in a unique way as  $\omega =\sum_{\nu =0}^{p^n -1}c_{\nu}w_{\nu}(dx)^{\otimes m},$  $c_{\nu}\in E$. 

We claim that $$(\sigma-1)^k \omega =0 \Rightarrow c_k =c_{k+1}=\cdots=c_{p^n -1}=0.$$
\textit{Proof of the claim}.
$(\sigma-1)^k \omega=(\sigma-1)^k \sum_{\nu =0}^{p^n -1}c_{\nu}w_{\nu}(dx)^{\otimes m}$. But $(\sigma-1)^k $ 
acts only on $w_{\nu}$ while leaving invariant the $c_{\nu}(dx)^{\otimes m}$. The $G$-action is given by
$(\sigma-1)^k w_{\nu}= 0,$ for all $k>\nu$. Indeed
$(\sigma-1)^{k+1} w_{k} =(\sigma -1)(\sigma-1)^k w_{k}=(\sigma -1)\prod_{\epsilon=1}^{n}\alpha_{\epsilon}^{(k)}!$
and $\alpha_{\epsilon}^{(k)}! \in K ,\;$ for all $\epsilon=1,\ldots,n$, so the product is fixed by the generator $\sigma$.
From that we have
\begin{eqnarray*}
(\sigma-1)^k \omega &=& (\sigma-1)^k \sum_{\nu =0}^{p^n -1}c_{\nu}w_{\nu}(dx)^{\otimes m}\\ 
%&=& \sum_{\nu =0}^{p^n -1}c_{\nu}(\sigma-1)^k w_{\nu}(dx)^{\otimes m}\\ 
&=& \sum_{\nu \geq k}^{p^n -1}c_{\nu}(\sigma-1)^k w_{\nu}(dx)^{\otimes m}\\ 
&=& c_k (dx)^{\otimes m} \prod_{\epsilon=1}^{n}\alpha_{\epsilon}^{(k)}! +\sum_{\nu \geq k+1}^{p^n -1}c_{\nu}(\sigma-1)^k w_{\nu}(dx)^{\otimes m}.
\end{eqnarray*} 
If the last equality is equal to zero, then we can see that 
\begin{eqnarray*}
(\sigma-1)^{k+1}\omega &=& 
(\sigma -1)c_{k}(dx)^{\otimes m}\prod_{\epsilon=1}^{n}\alpha_{\epsilon}^{(k)}! 
+\sum_{\nu \geq k+1}^{p^n -1}c_{\nu}(\sigma-1)^{k+1} w_{\nu}(dx)^{\otimes m} \\
&=&
\sum_{\nu \geq k+1}^{p^n -1}c_{\nu}(\sigma-1)^{k+1} w_{\nu}(dx)^{\otimes m}=0.
\end{eqnarray*}
We apply the above argument recursively and we finally get 
\[c_{p^n -1}=0.\] Now we write \[ \omega=\sum_{\nu =0}^{p^n -2}c_{\nu}(\sigma-1)^k w_{\nu}(dx)^{\otimes m},\] and we repeat the whole procedure to get  $ c_k =c_{k+1}=\cdots=c_{p^n -2}=0$ and prove the claim.

We now have an alternative expression for the quotients of $\Omega_{F}^{i}(m).$ Namely \linebreak
$\dim_K (\Omega_{F}^{k+1}(m)/\Omega_{F}^{k}(m))=$ 
\begin{equation}\label{dimpil}
=\dim_K \left\lbrace c_k \in E : \textrm{there is an } \omega \in \Omega_{F}(m), \textrm{  with  } \omega=
\sum_{\nu =0}^{k}c_{\nu} w_{\nu}(dx)^{\otimes m}\right\rbrace .
\end{equation}
%$\dim_K (\Omega_{F}^{k+1}(m)/\Omega_{F}^{k}(m))=$ 
%$$=\dim_K \left\lbrace c_k \in E : \textrm{there is an } \omega \in \Omega_{F}(m), \textrm{  with  } \omega=
%\sum_{\nu =0}^{k}c_{\nu} w_{\nu}(dx)^{\otimes m}\right\rbrace .$$
If $\omega =\sum_{\nu =0}^{k}c_{\nu}w_{\nu}(dx)^{\otimes m}$, then
$(\sigma-1)^k \omega= \prod_{\epsilon=1}^{n}\alpha_{\epsilon}^{(k)}! c_k (dx)^{\otimes m} \in \Omega_{F}(m),$
and from that we see that
$$ 0\leq v\left( i,n,\mu,\prod_{\epsilon=1}^{n}\alpha_{\epsilon}^{(k)}! c_k (dx)^{\otimes m}\right) =
v\left( i,n,\mu, c_k (dx)^{\otimes m}\right) ,$$
thus $c_{k} (dx)^{\otimes m}$ is a $G$-invariant $m$-holomorphic differential on $F$. This means that
\[
\mathrm{div}_{F} \left( c_k (dx)^{\otimes m}\right)=\mathrm{Con}_{F/E}\left(\mathrm{div}_E(c_k (dx)^{\otimes m})\right) + m \mathrm{Diff}(F/E)
\]
is an integral divisor. As the left side of the above equality has no poles,
the same will be true for the right side, so the only poles that we allow for $c_k (dx)^{\otimes m}$, are the ones that can be canceled out from the factor
$m\mathrm{Diff}(F/E)$. Hence $c_k (dx)^{\otimes m}$ can have poles only at ramified primes of $E$. As $\mathrm{div}(\omega)$
is an effective divisor, Lemma \ref{old2} requires that $v(i,n,\mu, c_k (dx)^{\otimes m} w_k) \geq 0,$ for all $\; i,\mu.$ So

\begin{equation}\label{bli}
v(i,n,\mu, w_k )+v(i,n,\mu,c_k (dx)^{\otimes m} ) \geq 0.
\end{equation}
We also have
\begin{eqnarray}\label{blii}
v(i,n,\mu, c_k (dx)^{\otimes m})&=&v(i,n,\mu,\mathrm{Con}_{F/E}\mathrm{div}_E(c_k (dx)^{\otimes m})) + v(i,n,\mu,m \mathrm{Diff}(F/E))\nonumber\\
&=& p^{e_i} \bar{v}_i (c_k (dx)^{\otimes m})+m \delta_i.
\end{eqnarray}
Combining Eq. (\ref{bli}) and (\ref{blii}), we obtain  $p^{e_i} \bar{v}_i (c_k (dx)^{\otimes m})+m \delta_i +v(i,n,\mu, w_k) \geq 0$, or
$$\bar{v}_i (c_k (dx)^{\otimes m}) \geq - \nu_{ik}(m).$$
The divisor of $E$ 
\begin{equation}\label{effective}
D=\sum_{i=1}^{s} \nu_{ik}(m)\bar{P}_i +\mathrm{div}_E(c_k ) +\mathrm{div}_E\left( (dx)^{\otimes m}\right)  \geq 0,
\end{equation}
is effective, which is equivalent to $c_k \in L\left(\sum_{i=1}^{s} \nu_{ik}(m)\bar{P}_i +\mathrm{div}_E\left( (dx)^{\otimes m}\right)  \right)$. Notice that the divisors $D$ and $\sum_{i=1}^{s} \nu_{ik}(m)\bar{P}_i +\mathrm{div}_E ( (dx)^{\otimes m})$ are linear equivalent. Thus they have the same degrees and $\ell(D)=\ell \left(\sum_{i=1}^{s} \nu_{ik}(m)\bar{P}_i +\mathrm{div}_E\left( (dx)^{\otimes m}\right)  \right)$. With that in mind, we will use the Riemann-Roch Theorem on the function field $E$ in order to 
compute
 $$\ell(D):=\dim_K L\left(\sum_{i=1}^{s} \nu_{ik}(m)\bar{P}_i +\textrm{div}_E\left( (dx)^{\otimes m}\right) \right).$$   We have
$\ell(mW)=\deg(mW)+1-g_E +\ell(W \setminus mW)$, where $W$ is a
 canonical divisor of $E$.
% ($W=(\omega$), for some $\omega \in \Omega_E $). 
It is well known that $\deg W=2g_E -2$, $\ell(W)=g_E$ and if $\deg(A)<0$ then $\ell( A)=0$, for every divisor $A$ of $E$.\\
We have the following cases:\\
\textit{Case 1}: $g_E \geq 2.$\\
 Hence $$\dim\Omega_E (m)=\ell(mW)=m(2g_E -2)+1-g_E +0$$ and $\ell(mW)=(2m-1)(g_E -1)$. Finally 
$$\ell\left(\sum_{i=1}^{s} \nu_{ik}(m)\bar{P}_i +\mathrm{div}_E\left((dx)^{\otimes m}\right) \right)=\deg\left(\sum_{i=1}^{s} \nu_{ik}(m)\bar{P}_i +\mathrm{div}_E\left((dx)^{\otimes m}\right) \right) +1 -g_E+ 0,$$ or equivalently $$\ell(D)=\Gamma_k (m)+(g_E -1)(2m-1).$$
\textit{Case 2}: $g_E =1.$\\
$\deg(W)=0$ and $\ell(mW)=1$, thus $$\ell(D)=\Gamma_k (m),$$ because
\begin{equation*}
\ell(D)=\Gamma_k (m) +m(2 g_E -2)+1 -g_E +\ell(W-mW-\Gamma_k (m)),
\end{equation*}
with the last term of the sum being zero because
\begin{equation*}
 \deg(W-mW-\Gamma_k (m))=(1-m)(2g_E -2)-\Gamma_k (m)=-\Gamma_k (m)<0.
\end{equation*}
\textit{Case 3}: $g_E =0.$\\
In that case $\deg(W)=-2<0$, thus $\ell(mW)=0$, for all $m\geq1$. Finally 
\begin{equation}\label{periorismos}
0\leq \ell(D)=\Gamma_k (m)-2m+1 ,
\end{equation}
because 
\begin{equation*}
\deg(D)=\Gamma_k (m)-2m\geq0, \textrm{ from Eq.(\ref{effective}). }
\end{equation*}
On the other hand as $c_k \in L\left(\sum_{i=1}^{s} \nu_{ik}(m)\bar{P}_i +\mathrm{div}_E\left((dx)^{\otimes m}\right)\right)$, from Eq. (\ref{dimpil}) we get, \textit{for all cases}, that
\begin{eqnarray}\label{inequality}
  \dim_K\left(\Omega_{F}^{k+1}(m)/\Omega_{F}^{k}(m)\right) \leq \ell(D)=\Gamma_k (m)+(g_E -1)(2m-1).
\end{eqnarray}

We want to show that inequality (\ref{inequality})  is actually an equality, so using Eq. (\ref{dks}), we are going to calculate the ${d_k} 's$. \\
We compute:
\begin{eqnarray}\label{from<to=}
(2m-1)(g_F -1)&=&\dim_K \Omega_F (m)=\sum_{k=0}^{p^{n}-1}\dim_K\left( \Omega_{F}^{k+1}(m)/\Omega_{F}^{k}(m) \right) \nonumber  \\
&\leq&\sum_{k=0}^{p^n -1}\Gamma_k (m) +\sum_{k=0}^{p^n -1}(g_E -1)(2m-1)\nonumber \\
&=&\sum_{k=0}^{p^n -1}\Gamma_k (m) +p^n (2m-1)(g_E -1).
\end{eqnarray}
If we show that $\dim_K \Omega_F (m)$ is equal to Eq. (\ref{from<to=}), then we have the desired equality. This is done using Riemann-Hurwitz Theorem:
$$2g_F -2 =\frac{[F:E]}{[K:K]}(2g_E -2) +\deg\mathrm{Diff}(F/E),$$
where $\mathrm{Diff}(F/E)=\sum_{i=1}^{s} \sum_{P(i,n,\mu )/\bar{P}_i}\delta_i \cdot P(i,n,\mu )$. Since $K$ is algebraically  closed, then $\deg P(i,n,\mu ) = f\left( P(i,n,\mu )/\bar{P}_i \right)=\left[ \mathcal{O}_{P(i,n,\mu )}/P(i,n,\mu ) :\mathcal{O}_{\bar{P}_i} /\bar{P}_i\right]=1,$ for all $i$. Moreover,  since the $F/E$ is  Galois the number of places of $F$ above $\bar{P}_i$ is $\mu$, where $[F:E]=e\left( P(i,n,\mu )/\bar{P}_i \right)\cdot f\left( P(i,n,\mu )/\bar{P}_i \right)\cdot \mu$, so $\mu=p^{n-e_i }$. We can now calculate
\begin{equation}\label{R-H}
(2m-1) (g_F -1 )=p^n (g_E -1)(2m-1) +\frac{1}{2}(2m-1)\sum_{i=1}^s p^{n-e_i } \delta_i .
\end{equation}
From Eq. (\ref{from<to=}) and (\ref{R-H}), it is enough to show that 
\begin{equation}\label{Differentalt}
\sum_{k=0}^{p^{n}-1}\nu_{ik} (m) =\frac{1}{2}(2m-1) p^{n-e_i }\delta_i , \  \textrm{ for all }i.
\end{equation}
\begin{remark}\label{degdif1}
Observe from Eq. (\ref{Differentalt}) that it is enough to show that
\begin{equation*}
\frac{2}{2m-1}\sum_{k=0}^{p^n-1}\Gamma_k(m) =\deg\mathrm{Diff}(F/E).
\end{equation*}

\end{remark}

We compute: 
\begin{eqnarray*}
\sum_{k=0}^{p^n -1} \nu_{ik}(m)&=&\sum_{k=0}^{p^n -1} \biggl\lfloor \frac{m\delta_i - \sum_{j=1}^{n}a_j ^{(k)}\Phi(i,j)p^{n-j}}{p^{e_i }}\biggr\rfloor\\
&=& p^{n-e_i }m\delta_i -\frac{1}{p^{e_i }}\sum_{k=0}^{p^{n}-1}\sum_{j=1}^{n}a_j ^{(k)}\Phi(i,j)p^{n-j}\\
&-&\sum_{k=0}^{p^n -1}\left\langle  \frac{m\delta_i - \sum_{j=1}^{n}a_j ^{(k)}\Phi(i,j)p^{n-j}}{p^{e_i }}\right\rangle.
\end{eqnarray*}
First we take care of  the second summation. Fix a $j$. As $k$ runs over $0,\ldots ,p^n -1,$  the elements $a_{j}^{(k)}$  take all the values from zero to $p-1$, $p^{n-1}$ times. Considering this we have 
\begin{eqnarray}
\frac{1}{p^{e_i }}\sum_{k=0}^{p^{n}-1}\sum_{j=1}^{n}a_j ^{(k)}\Phi(i,j)p^{n-j}&=& \frac{1}{p^{e_i }}\sum_{j=1}^{n}\Phi(i,j)p^{n-j}\sum_{k=0}^{p^{n}-1}a_j ^{(k)}\nonumber\\
&=&\frac{1}{p^{e_i }}\sum_{j=1}^{n}\Phi(i,j)p^{n-j} p^{n-1}\sum_{k=0}^{p-1}k\nonumber\\
&=&\frac{1}{2}p^{n-1}p(p-1)\frac{1}{p^{e_i }}\sum_{j=1}^{n}\Phi(i,j)p^{n-j}\nonumber\\
&=&\frac{p^{n-e_i }}{2} (p-1)\sum_{j=1}^{n}\Phi(i,j)p^{n-j}\nonumber\\
&=&\frac{p^{n-e_i }}{2}(\delta_i +1 -p^{e_i}), \label{sec.sum}
\end{eqnarray}
where the last equality came from Lemma \ref{different}.\\
Then we consider the fractional part. Observe that $\Phi(i,j)=0$, for all $j\leq n-e_i $ and $\Phi(i,j)$ is relatively prime to $p$, from the standard form hypothesis. 
We notice that as $a_{j}^{(k)}$ 's runs over $0,\ldots,p-1 $ for $j\geq n-e_{i} +1,$ the numbers  $\sum_{j=n-e_i +1}^{n} a_{j}^{(k)}p^{n-j}$, form a complete system $\mod p^{e_i }$. In the case where $r=0$, then the same numbers, for $j\geq 1,$ form a complete system $\mod p^{n }$. It is well known from elementary number theory that the same is true for $\sum_{j=n-e_i +1}^{n}a_{j}^{(k)}\Phi(i,j)p^{n-j}$, using the fact that $\textrm{g.c.d}(\Phi(i,j),p)=1$. Thus as $k$ runs over $0,\ldots,p^n -1,$ the numbers $m\delta_i - \sum_{j=1}^{n}a_{j} ^{(k)}\Phi(i,j)p^{n-j}$ run over a complete residue system $\mod p^{e_i }$ (in fact $z \pm \sum_{j=1}^{n}a_{j} ^{(k)}\Phi(i,j)p^{n-j}$ run through a complete  residue system for all $ z\in \mathbb{Z}$), $p^{n-e_i }$ times. We are ready to calculate

\begin{eqnarray}
\sum_{k=0}^{p^n -1} \left\langle\frac{m\delta_i - \sum_{j=1}^{n}a_j ^{(k)}\Phi(i,j)p^{n-j}}{p^{e_i }}\right\rangle&=&\frac{p^{n -e_i} }{p^{e_i }}\sum_{k=0 }^{p^{e_i } -1}k\nonumber\\
&=&\frac{p^{n -e_i} }{p^{e_i }}\frac{p^{e_i }(p^{e_i}-1)}{2}\nonumber\\
&=&\frac{1}{2}(p^{n}-p^{n-e_i })\label{frac.part}.
\end{eqnarray}
The final step is to combine the Equations (\ref{sec.sum}), (\ref{frac.part}) and the first summand, in order to obtain
\begin{eqnarray}\label{last}
\sum_{k=0}^{p^n -1} \nu_{ik}(m)&=&p^{n-e_i }m\delta_i -\frac{p^{n-e_i }}{2}(\delta_i +1 -p^{e_i} )-\frac{1}{2}(p^{n}-p^{n-e_i })\nonumber \\
&=&\frac{2 p^{n-e_i}m \delta_i -p^{n-e_i} \delta_i -p^{n-e_i}+p^{n} -p^{n}+p^{n-e_i}}{2}\nonumber \\
&=&\frac{1}{2}\left( p^{n-e_i} \delta_i (2m-1) \right),
\end{eqnarray}
We showed that inequality (\ref{from<to=}) is actually an equality. Using Eq. (\ref{dks}), to compute the $d_k$'s, Theorem follows.
\end{proof}

\begin{remark}[The case $m=1$, $r=0$]\label{m=1}
If $g_E\geq2$ then $\ell(D)=\Gamma_k (1) +(g_E -1)$, for all $0\leq k < p^n -1$ and $\ell(D)=(g_E -1) +1$ for $k=p^n -1$. If $g_E =1$,  then $\ell(D)=\Gamma_k (1)$, for $0\leq k \leq p^n-2$, while $\ell(D)=1$ for $k=p^n -1$. Finally if  $g_E =0$ then $\ell(D)=\Gamma_k (1) -1$, for $0\leq k \leq p^n-2$, while $\ell(D)=0$ for $k=p^n -1$. The extra cases that we do not consider in Eq.  (\ref{inequality}) are
\begin{eqnarray*}
\textrm{ for } g_E\geq2,\ d_{p^{n}}&=& \dim_K\left(\Omega_{F}^{p^n}(1)/\Omega_{F}^{p^n -1}(1)\right) \leq \ell(D)=g_E ,\\
\textrm{ for } g_E=1,\ d_{p^{n}}&= &\dim_K\left(\Omega_{F}^{p^n}(1)/\Omega_{F}^{p^n -1}(1)\right) \leq \ell(D)=1,\\
\textrm{ for } g_E=0,\ d_{p^{n}}&=& \dim_K\left(\Omega_{F}^{p^n}(1)/\Omega_{F}^{p^n -1}(1)\right) \leq \ell(D)=0.
\end{eqnarray*}

Observe that for the exceptional cases- the inequalities above-, we have that 
\begin{eqnarray*}
g_F &=&\dim_K \Omega_F (1)=\sum_{k=0}^{p^{n}-1}\dim_K\left( \Omega_{F}^{k+1}(m)/\Omega_{F}^{k}(m) \right) \nonumber  \\
&\leq&\sum_{k=0}^{p^n -1}\Gamma_k (1) +\bigl(p^n (g_E -1) +1\bigr).
\end{eqnarray*}
Thus, in order to  prove that  these inequalities are equalities, with the help of Eq.  (\ref{R-H}) we just have to show Eq. (\ref{Differentalt}) for $m=1$.

{This is a way to arrive at  \cite[Theorem 2]{vm}, when a place in $\mathbb{P}_E$ is totally ramified in $F$, 
or equivalently when there is not an unramified subextension of $F/E$  (i.e. $r=0$), in order to prove an analogous result for the $m=1$ case:

\begin{theorem}[Valentini-Madan]
For $m=1$ and when exists a place in $\mathbb{P}_E$ that is totally ramified in $F$,  the regular representation of $G$ occurs
 $d_{p^n}=g_E$ times in the representation of $G$ on $\Omega_F (1)$. For $k=1,\ldots, p^n -1$, the indecomposable representation
 of degree $k$ occurs $d_{k}=\Gamma_{k-1}(1) -\Gamma_{k}(1) +\beta$ times, where $\beta$ equals to $-1$, if $k=p^n-1$ and equals
 to  zero otherwise.
\end{theorem}
As Valentini and Madan observed for the case $m=1$, $\Gamma_k(1)=0$, if $\nu_{ik}(1)=0$ for all $i$'s, and that could happen if and only if
 $k\geq p^n-p^r$. However when $m\geq 2$ then $\Gamma_k(m)\gneqq 0$. Thus for the $m=1$ case one needs to distinguish cases 
on whether Boseck's invariants are zero or not, in order to derive  a result for the $K[G]$ module structure of $\Omega_F(1)$ which   will be
depended on $r$.  If we assume that $r\neq0$ in the case $m=1$ then one 
should use Tamagawa or Valentini results, \cite{Tamagawa:51} \cite{Val:82}, in  order to treat the unramified $E_r/E$ extension
 (see \cite[proof of Theorem 2 , case $k\geq p^n-p^r$]{vm}) .

This is the major difference between the exposition found there and the one followed here, which also shows  that for $m\gneq1$ 
Boseck invariants, and hence the conditions for holomorphicity that they compactly express, as well as the $K[G]$ module structure of $\Omega_F(m)$,
 do not depend on $r$ and thus on the existence of an unramified subextension, $E_r /E$, that appears  when no place is totally ramified in $F/E$.}
% We focus in the $m\geq2$ case, and treat the $m=1$ separately in Remarks, since that case . 

\end{remark}

\subsection{A New Basis for Holomorphic Differentials}\label{cyclic basis}
We now proceed to a basis construction for the holomorphic $m$-(poly)differentials $\Omega_{F}(m)$ when $E$ is rational. Without loss of generality, we will assume that $E=K(x)$ and $F/K(x)$ will be a cyclic extension of degree $p^n$. The main result of this subsection, looks like Lemma 5 of Madden \cite{Madden78}. The main difference is that there, he assumed that \textit{the infinite prime of $E$ ramifies} and takes two cases; case one stands for the $P(i,n,\mu)$'s that are lying over the ``finite'' primes of $E$ while case two, for $P(i,n,\mu)$'s that are lying over the  ``infinite'' prime of $E$, $p_{\infty}$.\\
Here we assume that \textit{the infinite prime of $E$ does not ramify} and give a single basis for the holomorphic differentials first ($m=1$), and then for holomorphic $m$-(poly)differentials, ($m>1$). 

The quantities that we defined at Eq. (\ref{niks}) play a major role for us. They are the quantities that Valentini and Madan defined at \cite[p.110]{vm} and the quantities that followed from the restrictions on $\nu$, at Lemma 5 of Madden, \cite{Madden78}. These are in fact, the same quantities that Boseck constructed in his work \cite{Boseck} in the late fifties, which led him to a similar basis for $\Omega_F (1)$, when $F/K(x)$ was an Artin-Schreier or a Kummer (of degree $q$, with g.c.d$(p,q)=1$) extension, of an algebraic function field. Let us have a closer analysis.

\[
\textrm{div}_{F}((dx)^{\otimes m})=m\mathrm{Diff}(F/K(x))-2m  (x)_{\infty},
\] 
or using Lemma \ref{different}, we have
\[m\sum_{i=1}^{s}p^{n-e_i}\bigr(  \sum_{j=1}^n (p-1) \Phi(i,j) p^{n-j}+(p^{e_i}-1) \bigl)\cdot P(i,n,\mu) -2m\textrm{Con}_{F/K(x)}(p_{\infty})= \]
\begin{equation*} %\label{divdx}
=m\sum_{i=1}^{s}p^{n-e_i}\delta_{i}\cdot P(i,n,\mu)-2m\textrm{Con}_{F/K(x)}(p_{\infty}).
\end{equation*}
Using Lemma \ref{different} and observing that $v(i,n,\mu,w_k)=v(i,n,\mu',w_k)$, we can write: %\footnote{here i had forgoten $p^{n-e_i}$ in the following sum }
\begin{equation*}
\textrm{div}_{F}(w_k )=P_{w_k} -\sum_{i=1}^{s}p^{n-e_i}\sum_{j=1}^n a_j^{(k)} \Phi(i,j) p^{n-j}\cdot P(i,n,\mu),
\end{equation*}
where $P_{w_k}$ is an effective divisor of $F$, prime to $ P(i,n,\mu)$, for all $i$'s. Then

\begin{equation*}
\textrm{div}_{F}(w_k (dx)^{\otimes m})=P_{w_k} +\sum_{i=1}^{s}p^{n-e_i}\left(m\delta_{i}-\sum_{j=1}^n a_j^{(k)} \Phi(i,j) p^{n-j}\right) \cdot P(i,n,\mu)-2m\textrm{Con}_{F/K(x)}(p_{\infty}).
\end{equation*}

We analyze $p^{n-e_i}\bigl(m\delta_{i}-\sum_{j=1}^n a_j^{(k)} \Phi(i,j) p^{n-j}\bigr)$ $\mod p^{n}$:
\begin{equation*}
p^{n-e_i}\bigl(m\delta_{i}-\sum_{j=1}^n a_j^{(k)} \Phi(i,j) p^{n-j}\bigr)=\nu_{ik}(m) \cdot p^{n}+\rho_{i}^{(k,m)}, \textrm{ with } 0\leq \rho_{i}^{(k,m)}\leq p^{n} -1.
\end{equation*}
Notice that we are consistent with the definition of $\nu_{ik}(m)$'s given in Eq. (\ref{niks}). We denote the ramified places of E, $\bar{P_{i}}, i=1,\ldots ,s $ with $\bar{P_{i}}(x),$  (since $E$ is rational then every ramified place of $E, \bar{P_{i}}$, corresponds to an irreducible polynomial, $\bar{P_{i}}(x)$ )  and let
\begin{eqnarray*}
 g_k (x)=\prod_{i=1}^{s} \bar{P_{i}}(x)^{\nu_{ik}(m)}.
\end{eqnarray*}
Then the quantity $\Gamma_k (m)$ of Eq. (\ref{niks}) is naturally defined as the degree of $g_k(x)$. This is the Boseck invariant for this case. We then have
\begin{equation}\label{divbases}
\textrm{div}_{F}(w_k [g_{k}(x)]^{-1} (dx)^{\otimes m})=\sum_{i=1}^{s}\rho_{i}^{(k,m)}\cdot P(i,n,\mu)+\left(\Gamma_k (m)-2m\right)\textrm{Con}_{F/K(x)}(p_{\infty})+P_{w_k} .
\end{equation}
Observe that if 
\begin{equation}\label{Gk basis}
\Gamma_k (m)\geq2m,
\end{equation}
then, the divisor of Eq. (\ref{divbases}) is effective, so for $k=0\ldots,p^n -1 ,\    \Gamma_k (m)\geq2m$ and for $0\leq \nu \leq\Gamma_k (m)-2m$, the  $x^\nu w_k [g_{k}(x)]^{-1} (dx)^{\otimes m}$ is a holomorphic $m$-(poly)differential of $F$.

Recall that the set $w_k=y_1^{a_1^{(k)}}y_2^{a_2^{(k)}}\cdots y_n^{a_n^{(k)}}$, $0\leq k \leq p^n-1$ is an  $E$-basis of $F$. The main result of this subsection is:
\begin{lemma}\label{basis generalized A-S}
The set 
$\langle x^\nu w_k [g_{k}(x)]^{-1} (dx)^{\otimes m} : 0\leq k  \leq p^n -1 +\beta_m, 0\leq \nu \leq \Gamma_k(m) -2m , m\geq1 \rangle$, is a $K$-basis for the space of the $m$-(poly)differentials $\Omega_F (m)$, with $\beta_1 =-1$, if $m=1$ and $\beta_m=0$ if $m\geq2$.
%\footnote{i corrected the basis in the case where m=1} 
\end{lemma}
\begin{proof}
 First remember that
\begin{equation*}
\nu_{ik}(m):=\biggl\lfloor\frac{m\delta_i -\sum_{j=1}^n a_j^{(k)} \Phi(i,j) p^{n-j}}{p^{e_i}}\biggr\rfloor.
\end{equation*}
We now observe that:\\
Since the base field is rational,  $r=0$, because every finite separable extension of the rational function field should be ramified. Thus, for $m=1$, the condition  $a_j^{(k)}=p-1$, for all $j=1,\ldots n$, is equivalent to $k=p^n -1$ and $\Gamma_{p^n -1} (1)=0$, while\\
for $m>1$, we always have that $\Gamma_k (m)\neq 0$ for every $k$.\\ 
It is enough to prove the equalities below:
\begin{enumerate}
\item[(i)] For $m=1$, we must prove that $\sum_{k=0}^{p^n -2}\bigl( \Gamma_k (1)-1\bigr)=\dim_K \Omega_F (1) =g_F$,
\item[(ii)]  For $m\geq2$, we must see that $\sum_{k=0}^{p^n -1}\bigl( \Gamma_k (m)-2m+1\bigr)=\dim_K \Omega_F (m)=(2m-1)(g_F -1)$. This number is well defined since $g_F \geq 2$ from our hypothesis.
\end{enumerate}
For case (i),
%\footnote{i throwed Maddan's reference and put mine}
 we see from Eq. (\ref{last}) for $m=1$ that this number has been computed :
$$\sum_{k=0}^{p^n -1} \nu_{ik}(m)=\frac{1}{2}\left( p^{n-e_i} \delta_i  \right).$$
So
\begin{equation*}
\sum_{k=0}^{p^n -2}\bigl( \Gamma_k (1)-1\bigr)=\frac{1}{2}\sum_{i=1}^{s}\left( p^{n-e_i} \delta_i  \right)-p^n+1=g_F .
%\sum_{k=0}^{p^n -2}\bigl( \Gamma_k (1)-1\bigr)=\frac{1}{2}\sum_{i=1}^{s}p^{n-e_i}[(p-1)  \sum_{j=1}^n \Phi(i,j) p^{n-j}+p^{e_i}-1]-p^n+1=g_F .
\end{equation*}
Where the last equality comes from Riemann-Hurwitz Theorem, a version of that can be found in Eq. (\ref{R-H}).  For $m \geq 2$  Eq. (\ref{last}), gives that
\begin{eqnarray*}
\sum_{k=0}^{p^n -1}\nu_{ik}(m)&=&\frac{1}{2}\left( p^{n-e_i} \delta_i (2m-1) \right), \textrm{ thus }\\
\sum_{k=0}^{p^n -1}\bigl( \Gamma_k (m)-2m+1\bigr)&=&\sum_{i=1}^{s}\bigl(\frac{1}{2} p^{n-e_i} \delta_i (2m-1)\bigr) -p^n (2m-1)\\
&=&(2m-1)\left(\frac{1}{2}\sum_{i=1}^{s} p^{n-e_i} \delta_i - p^n\right)\\
&\stackrel{\textrm{Eq. } (\ref{R-H})}{=}&(g_F -1)(2m-1).
\end{eqnarray*}
\end{proof}
The following Remarks are here to confirm the correctness of our arguments, given in the proof of the Theorem \ref{basic thrm}, in subsection \ref{Galoismodulecyclic}.
\begin{remark}
Observe that the restriction for $\Gamma_k (m)$, given in Eq. (\ref{periorismos}), when the genus of $E$ was zero, is the same restriction that results from our basis in order to ensure that an $m$-(poly)differential is holomorphic (see Eq. (\ref{Gk basis})).
\end{remark}
\begin{remark}\label{2proof}
A second, simpler proof for Theorem \ref{basic thrm} can be given when $g_E =0$, using the basis of Lemma \ref{basis generalized A-S}. Indeed having Eq. (\ref{dimpil})  instead of using Riemann-Roch Theorem to count  the $c_k$'s, we use our basis. %The $c_k \in E$ such that $\omega=\sum_{\nu=0}^{k}c_{\nu}w_{\nu}(dx)^{\otimes m}$ 
From Eq. (\ref{dimpil}) the needed differentials are of the form $\omega= x^k [g_k (x)]^{-1}w_{k}(dx)^{\otimes m}$  with $c_k :=x^k [g_k (x)]^{-1}$ and there are $\Gamma_k (m) -2m +1$ of them, for every $k=1,\ldots, p^n -1$. So for  $k=1,\ldots, p^n-1$, with the help of Eq. (\ref{dks}), we get that $d_k = \Gamma_{k-1} (m) -2m +1-(\Gamma_{k}(m)-2m +1)=\Gamma_{k-1} (m)-\Gamma_{k} (m)$, while for $k=p^n$ we have that $d_{p^n}=\Gamma_{p^n-1}(m) -2m+1$.
\end{remark}
%\footnote{i added this paragraph plus the }

\subsection{A classical theorem of Hurwitz} \label{HuClass}
We will study  now  how Boseck invariants behave when $F/E$ is a cyclic tame ramified extension of function fields. Let $F/E$ is be a cyclic Kummer extension of degree $n$, with  $\textrm{g.c.d}(n,p)=1$. Choose a primitive $n$th root of unity, $\zeta$ of $K$. Let $F=E(y)$ and 
\begin{equation}\label{kumeq}
y^n =u, u\in E.
\end{equation}
Let $\bar{P}_i \in \mathbb{P}_E $, with $1\leq i \leq r$ , be the ramified places of $E$ in $F$ and $P_i \in \mathbb{P}_F$ the places above $\bar{P}_i$. We can assume that if $\sigma$ is a generator of $G$ then $\sigma(y)=\zeta y$ and also that $0\lneqq v_{\bar{P}_i}(u) \lneqq n$. Using Kummer theory \cite[p.110, Proposition III.7.3]{StiBo} we have that $e_i=e\bigl(P_i/\bar{P}_i \bigr)=\frac{n}{\textrm{g.c.d.}(n, v_{\bar{P}_i}(u))}$ and $\delta_i =e_i -1$. Finally, set $\Phi(i)=v_{P_i}(y) =\frac{e_i v_{\bar{P}_i}(u)}{n}$.

If we define $\Omega_{F}^{k}(1) =\{\omega \in \Omega_F(1) \ | \sigma(\omega)=\zeta^k \omega\},$ for $0\leq k \leq n-1$, then
\[
\Omega_F(1)=\bigoplus_{k=0}^{n-1}\Omega_{F}^{k}(1),
\]
and the corresponding $d_k$'s for this case are $d_k =\dim_K\Omega_{F}^{k}(1)$. It is also well known (see at \cite[p.272, V.2.2, Corollary]{Farkas-Kra}) that $d_0 =\dim_K\Omega_{F}^{G}(1)=g_E$. Notice that this is false  in positive characteristic in the 
case of wild ramification. 

In order to find the $K[G]$ module structure of $\Omega_F(1)$, we should compute once more the $d_k$'s. 
For that reason we should find the Boseck invariant for this case, i.e.  consider the corresponding rational extension. If $E=K(x)$,  then Boseck 
\cite[p. 50, Satz 16]{Boseck}, proved that (for $m=1$) :
\begin{proposition}[Boseck]\label{basisKummer}
The set $\langle x^\nu y^{-k} g_{k}(x)dx : 0\lneqq k \leq n-1  , 0\leq \nu \leq \Gamma_k( 1) -2 \rangle$ is a $K$-basis for the space of  $\Omega_F (1)$, where $g_{k}(x)=\prod_{i=1}^{r}\bar{P}_i (x)^{\lfloor \frac{k\Phi(i)}{e_i}\rfloor}$ and
\begin{equation*} \label{Gamma_k for the tame case}
 \Gamma_k( 1)=\sum_{i=1}^r\left\langle \frac{k\Phi(i)}{e_i}\right\rangle.
\end{equation*}
\end{proposition}
\begin{remark}\label{degdifKummer}
One crucial step in the proof of Proposition \ref{basisKummer} is to show that
\begin{equation*}
2\sum_{k=1}^{n-1}\Gamma_k( 1)=\deg\mathrm{Diff}(F/E).
\end{equation*}

\end{remark}

We have $k$ distinct irreducible representations of degree $1$. If $g_E=0$ and $k\neq 0$, we may count the differentials in Proposition \ref{basisKummer} in order to find the $d_k$'s. We see that $\omega \in \Omega_{F}^{k}(1)$ if $\omega=x^\nu y^{n-k} g_{n-k}(x)dx$. Their  number equals to 
\begin{eqnarray*}
\Gamma_{n-k}( 1) -1&=&\sum_{i=1}^r\left\langle \frac{n-k\Phi(i)}{e_i}\right\rangle-1\\
&=&\sum_{i=1}^r\left\langle \frac{-k\Phi(i)}{e_i}\right\rangle-1.
\end{eqnarray*}
So the $k$th representation occurs $d_{n-k}=\Gamma_{n-k}( 1) -1$ times in the representation of $G$ in $\Omega_F(1)$.
When $k=0$ we know that occurs $d_0=g_E=0$ times. 

{Notice that, in general, the $\Gamma_k (m)$'s depend on the different exponent and  the evaluation of basis elements, thus the genus of the base field  $E$,
  $g_E$ does not affect them.} We can  claim now that the same result will be true when $g_E\geq0$. Indeed this is the Hurwitz theorem (compare also to \cite[Theorem 2]{vm}):
\begin{theorem}[Hurwitz]\label{cyclic prime to p}
% \footnote{here the regular representation occurs $d_0=g_E=d_n$ times and each $kth$ irreducible representation is contained one time (``with multiplicity one'') in the regular representation, This is the reason we add the genus in the above formula}
For $k=0,\ldots, n -1$, we have $n$ distinct  irreducible representations of degree $1$. The $k$th representation occurs $d_{n-k}:=\Gamma_{n-k}(1) -1+g_E$ times in the representation of $G$ in $\Omega_F(1)$, when $k\neq 0$ and $g_E$ times when $k=0$.
\end{theorem}
% One should expect similar results when $m>1$.
%
%
%
\section{The Elementary Abelian Case}\label{elab}

Let $K$ be an algebraically closed field of characteristic $p$ and consider an elementary abelian 
Galois extension $F/K(x)$ of the rational function field $E=K(x)$. Set $G=\mathrm{Gal}(F/K(x))$ and assume that $|G|=p^n$. We also assume  that every place of $K(x)$ that is ramified in the above extension is totally ramified, i.e, $F=K(x, y)$, where $y$ satisfies the equation:
\begin{eqnarray}\label{equation e.abelian}
y^{q}-y =\frac{g(x)}{(x-a_1 )^{\Phi(1) }\cdots(x-a_{r})^{\Phi(r) }},
\end{eqnarray}
where $q=p^n $, $n\geq 1, \  g(x)\in K[ x ], \textrm{ and } \deg g(x)\leq \sum_{i=1}^{r}\Phi(i) =M,$ $g(a_i )\neq 0$. Thus, $p_{\infty}$, the infinite place of $K(x)$, is assumed to be unramified in $F$. Finally all the $\Phi(i) $'s are relatively prime to $p$ for all $i=1\ldots, r$. Let $p_i \in \mathbb{P}_{K(x)},\; i =1, \ldots, r$, be the rational places of $K(x)$ that are totally ramified in $F$, corresponding to  $(x-a_i ), \ a_i \in K$. 

First we fix some notation, similar to the notation given in subsection \ref{cyclic basis}:\\
It is well known (see for example Stichtenoth, \cite[ Prop. III.7.10]{StiBo}) that $\mathrm{Con}_{F/K(x)}(p_{i})=p^n P_{i}$, for $P_i \in \mathbb{P}_F$ above $p_i$ and $\mathrm{Diff}(F/K(x))= \sum_{i=1}^{r}(p^n -1)(\Phi(i) +1) P_i $, for all $i=1,\ldots,r$. We also know that for any $\beta \in \mathbb{F}_q$ the element $\sigma_{\beta}\in G$  acts  on the generator of the extension $y$ as: $\sigma_{\beta}(y)=y+\beta$. For an $m$-(poly)differential of $F$, $(dx)^{\otimes m}$ we will also have
\begin{eqnarray}\label{divdx^m}
 \textrm{div}_{F}((dx)^{\otimes m})&=&m\mathrm{Diff}(F/K(x))-2m  (x)_{\infty}\nonumber\\
&=&m\sum_{i=1}^{r}(p^n -1)(\Phi(i) +1) P_i -2m\textrm{Con}_{F/K(x)}(p_{\infty}),
\end{eqnarray}
where $(x)_{\infty}$ is the pole divisor of $x$. Taking degrees to the above divisors (or alternatively using the Riemann-Hurwitz formula), we  have 
\begin{eqnarray}\label{genus}
g_F =\frac{p^n -1}{2}\left(-2+\sum_{i=1}^{r}(\Phi(i) +1)\right).
\end{eqnarray}
Set $\textrm{div}_F(y)=P_y -\sum_{i=1}^{r}\Phi(i) P_i $, for an effective divisor $P_y$ of $ F$, prime to $P_i$, for $i=1,\ldots,r$. Then for $k=0,\ldots, p^n -1$ we compute $\textrm{div}_{F}  (y^{k} (dx)^{\otimes m})$, which equals to
\begin{eqnarray}
&=&\textrm{div}_F (y^{k}) +\textrm{div}_F ((dx)^{\otimes m})\nonumber \\ 
&=&k P_{y} -k\sum_{i=1}^{r} \Phi(i) P_i +m\sum_{i=1}^{r}(p^n -1)(\Phi(i) +1) P_i -2m\textrm{Con}_{F/K(x)}(p_{\infty})\nonumber\\ 
&=&\sum_{i=1}^{r}\biggr(\bigr((p^n -1)m -k\bigl)\Phi(i) +m(p^n -1)\biggl) P_i +k P_{y} -2m\textrm{Con}_{F/K(x)}(p_{\infty}).
\end{eqnarray}
Then we consider $\bigr((p^n -1)m -k\bigl)\Phi(i) +m(p^n -1)$ modulo $p^n$ 
$$\bigr((p^n -1)m -k\bigl)\Phi(i) +m(p^n -1)=\nu_{ik}( m)\cdot p^n +\rho_{i}^{(k, m)},$$
where $0\leq\rho_{i}^{(k, m)}\leq p^n -1$, is the remainder of the division, so
\begin{equation}\label{mi}
 \nu_{ik}( m)=\biggr \lfloor\frac{m(p^n -1 )(\Phi(i) +1)-k \Phi(i) }{p^n}\biggl\rfloor.
\end{equation}
We also take 
\begin{equation}\label{gkx}
g_{k}(x):=\prod_{i=1}^{r}(x-a_i )^{\nu_{ik}( m) },
\end{equation}
and let $\Gamma_k (m)=\sum_{i=1}^{r} \nu_{ik}(m),$ be the Boseck invariant for this case. 

Observe from Eq. (\ref{mi}) that for $m=1$, if $\Gamma_k(1)=0$, then $\nu_{ik}( 1)=0$ 
for all $i=1,\ldots,r$  and that could happen only when $k =p^n -1$ since $\nu_{ik}( 1) \geq 0$. 

For $m>1$, this is not the case because  $\nu_{ik} ( m)>0,$ for all $k$'s. 
We finally observe that $\textrm{div}_{F}[g_{k}(x)]^{-1}y^{k} (dx)^{\otimes m}$ is an effective divisor if 
\begin{equation}\label{holomdif}
 \Gamma_k (m)=\sum_{i=1}^{r} \nu_{ik}(m)\geq 2m,
\end{equation}
because
\begin{equation*}
\textrm{div}_{F}  (y^{k}[g_{k}(x)]^{-1} (dx)^{\otimes m})=\sum_{i=1}^{r}\rho_{i}^{(k, m)} P_i +k P_{y} + \bigr(\sum_{i=1}^{r} \nu_{ik}(m)-2m\bigl)\textrm{Con}_{F/K(x)}(p_{\infty}).
\end{equation*}
So when inequality (\ref{holomdif}) is fulfilled, then $x^\nu y^k [g_{k}(x)]^{-1} (dx)^{\otimes m}$ is a holomorphic $m$-(poly)differential, for $0\leq k \leq p^n -1 , 0\leq \nu \leq\Gamma_k (m) -2m$ and $m\geq1 $.
\begin{remark}\label{unifying remark}
We kept the notation of this section as close is possible, to the notation used in Section 2. Generally the quantities $\nu_{ik}(m)$, for a $p$-extension with $p|\textrm{char}K$, are equal to
\begin{equation*}
\left\lfloor\frac{m\delta_i + \{ \textrm{evaluation of the kth  E basis element of F  by a normalized valuation of F }\} }{p^{e_i}}\right\rfloor,
%\footnote{i put kth in place of an``an'' plus i changed the Remark a little}
\end{equation*}
where i runs the ramified places of $E$, $\delta_i$ is the different exponent of the extension and $e_i$ the ramification indices of the ramified primes of $E$ in $F$. The basis element is evaluated by a (normalized) valuation determined by a place of $F$ above a ramified place of $E$. Here we have total ramification, so $e_i=n$. In the case where $F/E$ was a cyclic extension, the term in the brackets is nothing else than $v(i,n,\mu,w_k)=-\sum_{j=1}^n a_j^{(k)} \Phi(i,j) p^{n-j}$ given in Lemma \ref{different}. Here if we take the standard $E$ basis for $F$, namely, $\{w_k:= y^k |\ 0\leq k \leq p^n -1, \}_E$ , we will then have that $v_{P_i}(w_k)=kv_{P_i}(y)$. Letting the right hand of Eq. (\ref{equation e.abelian}) be equal to an $u\in E$, then from the strict triangle inequality \cite[p.5, Lemma I.1.10]{StiBo} we have that $v_{p_i}(u)=\min \{v_{p_i}(y^q), v_{p_i}(y)\}=qv_{p_i}(y)=v_{P_i}(y)$. Letting now $-\Phi(i)=v_{p_i}(u)$ the Remark follows. Finally observe that $\Phi(i,j)=\Phi(i)$ because we haven't got a cyclic extension, thus we cannot get an analogous field tower like we had in Eq. (\ref{Artin-inter}). That, explains the independence from $j$.
\end{remark}

\subsection{Basis Construction}

Following Boseck (see \cite{Boseck}, Satz 15), we now prove the analogue of Lemma \ref{basis generalized A-S}.
\begin{proposition}\label{Satz15}
The set $\Sigma_m :=\langle x^\nu y^k [g_{k}(x)]^{-1} (dx)^{\otimes m} : 0\leq k \leq p^n -1 +\beta_m , 0\leq \nu \leq \Gamma_k( m) -2m , m\geq1 \rangle$ is a $K$-basis for the space of the m-(poly)differentials $\Omega_F (m)$, with $\beta_1 =-1$, if $m=1$ and $\beta_m=0$ if $m\geq2$.
\end{proposition}
\begin{proof}
For $m=1$, we have nothing to prove since this is Theorem 2 of Garcia  \cite{garciaelab}. For $m> 1 $ it is enough to show that $\sum_{k=0}^{p^n -1}(\Gamma_k( m) -2m+1)=\dim_K \Omega_F (m).$ 

First, we will compute $\sum_{k=0}^{p^n -1} \nu_{ik}(m)$. This is equal to 
\begin{eqnarray*}
&=&\sum_{k=0}^{p^n -1}\frac{m(p^n -1 )(\Phi(i) +1)-k \Phi(i)}{p^n}- \sum_{k=0}^{p^n -1}\left\langle \frac{m(p^n -1 )(\Phi(i) +1)-k \Phi(i)}{p^n}\right\rangle \\
&=&\sum_{k=0}^{p^n -1}\frac{m(p^n -1 )(\Phi(i) +1)}{p^n}-\sum_{k=0}^{p^n -1}\frac{k \Phi(i) }{p^n}-\sum_{i=1}^{p^n-1}\left\langle\frac{m(p^n -1 )(\Phi(i) +1)-k \Phi(i)}{p^n}\right\rangle  \\
&=&m(p^n -1 )(\Phi(i) +1)-\frac{\Phi(i) (p^n -1)}{2}-\sum_{k=0}^{p^n -1}\left\langle\frac{m(p^n -1 )(\Phi(i) +1)-k \Phi(i)}{p^n}\right\rangle.
\end{eqnarray*}
As $k$ runs  a complete residue system $\mod p^n$ and $\textrm{g.c.d.} (\Phi(i) , p)=1$, the same is true for $m(p^n -1 )(\Phi(i) +1)-k \Phi(i) ,$ so
\begin{eqnarray}\label{profinal}
\sum_{k=0}^{p^n -1} \nu_{ik}(m)&=&m(p^n -1 )(\Phi(i) +1)-\frac{\Phi(i) (p^n -1)}{2}-\frac{1}{p^n}\sum_{k=0}^{p^n -1}k \nonumber\\
&=&m(p^n -1 )(\Phi(i) +1)-\frac{\Phi(i) (p^n -1)}{2}-\frac{p^n -1}{2}\nonumber \\
&=&m(p^n -1 )(\Phi(i) +1)-\frac{ (p^n -1)}{2}(\Phi(i) +1)\nonumber\\
&=&(2m-1)(\Phi(i) +1)\frac{(p^n -1)}{2}.
\end{eqnarray}
\begin{remark}\label{degdif3}
Observe, for one more time that from Eq. (\ref{profinal}) we have that
\begin{equation*}
 \frac{2}{2m-1}\sum_{k=0}^{p^n-1}\Gamma_k(m)=\deg\mathrm{Diff}(F/E).
\end{equation*}
\end{remark}
Using Eq. (\ref{profinal}) and (\ref{genus}), we will then have that 
\begin{eqnarray}\label{final}
 \sum_{k=0}^{p^n -1}\Gamma_k (m)&=& \sum_{k=0}^{p^n -1} \sum_{i=1}^{r}\nu_{ik}(m)\nonumber\\
&=&(2m-1)(g_F+p^n -1).
\end{eqnarray}
Finally from Eq. (\ref{final}) we see that 
\begin{eqnarray*}
\sum_{k=0}^{p^n -1}(\Gamma_k (m)-2m+1)&=&\sum_{k=0}^{p^n -1}\Gamma_k (m)-2m p^n +p^n \\
&=&(2m-1)(g_F+p^n -1)-p^n (2m-1)\\
&=&(2m-1)(g_F -1).
\end{eqnarray*}
\end{proof}
\begin{remark}
If $F/K(x)$ is an extension of degree $p$ ( i.e. an Artin Schreier)  then  Boseck, \cite{Boseck}  proved that the set $\Sigma_1$  is indeed a basis for $\Omega_F(1)$.  With exactly the same arguments used here, one can show that the set defined in Proposition \ref{Satz15}, is still a basis of $\Omega_F (m)$ when $F/K(x)$ is a cyclic Artin Schreier extension of degree $p$.
\end{remark}

\begin{remark}\label{wp}
A basis for the holomorphic $m$-(poly)differentials is closely connected with the computation of the $m$-Weierstrass points of $F$. The case where $m=1$, see Garcia and Boseck ( \cite{garciaelab}, \cite{garciaa-s} and \cite{Boseck}), leads to the computation of (classical) Weierstrass points. One can now follow their ideas for the case $m>1$.
\end{remark}

\subsection{Galois module structure of  \boldmath$ \Omega_F (m)$}\label{elabgaloisstructure}

In this subsection we will determine the Galois module structure of the space $\Omega_F (m)$, of $m$-(poly)differentials using Proposition \ref{Satz15}. Following closely the ideas of Calder\'{o}n, Salvador and Madan (\cite{csm}, Theorem 1), who considered the case $m=1$, we will generalize their result for $m>1$.

Let $\theta_0 ,\ldots, \theta_{p^n -1} \in K$ and for $j=1, \ldots, p^n $, let $W_j =\langle \theta_0 ,\ldots, \theta_{j-1}\rangle_{K}$. Let also $G$ acting on these $\theta_i$'s with the action described as follows:
\begin{equation*}
\sigma_{\alpha}(\theta_i )=\sum_{l=0}^{i}\left(\begin{array}{clr}i\\ l \end{array}\right)\alpha^{i-l}\theta_l ,\textrm{ for }\  0\leq i \leq k.
\end{equation*}
In the near future, we will interpret these $\theta_i$'s as sums of specific $m$-holomorphic (poly)differentials (the anxious reader should jump to Eq. (\ref{theta's}) ). 
\begin{theorem}\label{fff}
Let $F/K(x)$ as above. The $K[G]$ module $W_j$, for $j=1, \ldots, p^n$ is indecomposable, and
$$\Omega_F (m)\simeq \bigoplus_{j=1}^{p^n } W_{j}^{d_j},$$
where $d_{p^n }$ equals to $\Gamma_{p^n -1}( m)-2m +1$ and is the number of times that the regular representation of $G$ occurs in $\Omega_F (m)$, while the indecomposable representation $W_j$ occurs $d_j =\Gamma_{j-1}(m)-\Gamma_{j}(m)$ 
times, for all $j=1, \ldots, p^n -1,$ and for all $m> 1$.
\end{theorem}
\begin{proof}
Let $\omega_{k, \nu}^{m}$ be an arbitrary basis element of $ \Omega_{F}(m)$ 
and $\sigma_{\alpha}\in G$ for an $\alpha \in \mathbb{F}_{q}$. Then 

\begin{equation}\label{action}
\sigma_{\alpha}(\omega_{k, \nu}^{m})=x^\nu (y+\alpha)^k [g_{k}(x)]^{-1} (dx)^{\otimes m}=\sum_{i=0}^{k} \left(\begin{array}{clr}k\\ i \end{array}\right) \alpha^{k-i}y^i x^{\nu}[g_{k}(x)]^{-1} (dx)^{\otimes m}.
\end{equation}
For $0\leq i\leq k,$ let $h_{i,k}(x)=\prod_{j=1}^{r}(x-a_j )^{\nu_{ji} (m)-\nu_{jk}(m)}.$ As $0\leq i\leq k,$ we have that $\nu_{ji} (m)\geq \nu_{jk}(m)$. Thus
\begin{eqnarray}\label{polynomial}
K[x]\ni h_{i,k}(x)=\sum_{e=0}^{n(i,k, m)}b_{e}^{(i, k, m)}x^e ,
\end{eqnarray}
with $n(i,k, m)= \Gamma_i (m)- \Gamma_k (m)$ and $b_{n(i,k, m)}^{(i, k, m)}=1.$ From Equations (\ref{action}), (\ref{polynomial}) and (\ref{gkx}) we obtain that
\begin{eqnarray*}
\sigma_{\alpha}(\omega_{k, \nu}^{m})&=&\sum_{i=0}^{k} \left(\begin{array}{clr}k\\ i \end{array}\right) \alpha^{k-i}y^i x^{\nu}[g_{i}(x)]^{-1} h_{i,k}(x) (dx)^{\otimes m}\\
&=&\sum_{i=0}^{k}\sum_{e=0}^{n(i,k, m)} \left(\begin{array}{clr}k\\ i \end{array}\right) \alpha^{k-i}y^i x^{\nu+e}[g_{i}(x)]^{-1} b_{e}^{(i, k, m)} (dx)^{\otimes m}.
\end{eqnarray*}
Since $0 \leq e\leq  \Gamma_i(m)- \Gamma_k (m)$, we have that $0\leq \nu+e \leq   \Gamma_i (m) -2m$. This means that $y^i x^{\nu+e}[g_{i}(x)]^{-1}(dx)^{\otimes m}$ is a basis element, namely $\omega_{i, \nu+e}^{m}.$ Thus, we have that the $G$ action on the basis of $\Omega_F (m)$, is given by
\begin{eqnarray}\label{gaction on the basis}
\sigma_{\alpha}(\omega_{k, \nu}^{m})&=&\sum_{i=0}^{k}\sum_{e=0}^{n(i,k, m)} \left(\begin{array}{clr}k\\ i \end{array}\right) \alpha^{k-i} b_{e}^{(i, k, m)}\omega_{i, \nu+e}^{m}\nonumber \\ 
&=&\sum_{i=0}^{k}\left(\begin{array}{clr}k\\ i \end{array}\right) \alpha^{k-i}\left \{\sum_{e=0}^{n(i,k, m)} b_{e}^{(i, k, m)}\omega_{i, \nu+e}^{m}\right\}.
\end{eqnarray}
Observe that the coefficient of $\omega_{k,\nu}^m$
in the right side of Eq. (\ref{gaction on the basis}) is 1.\\
Let $M_{k,\nu}^{m}$, be the $K[G]$-module, generated by 
\[\{ \omega_{i, \nu+e}^{m} | 0\leq i\leq k,\ 0\leq e\leq n(i,k,m)\}.\]
The following cases play an important role in the decomposition of $M_{k,\nu}^{m}$, which will follow:
%Observe that $e$ takes the value zero twice:
\begin{itemize}
\item The condition $i=k$, implies $e=0$ and the differentials of $M_{k,\nu}^{m}$ which satisfy that condition, are of the form $\{\omega_{k,\nu}^m\}$, while
\item The conditions $e=0$ and $i\neq k$ are satisfied by the differentials of $M_{k,\nu}^{m}$, of the form $\{\omega_{i,\nu}^m| \ 0\leq i\lneqq k\}$.
\end{itemize}
With these in mind, we have that
\begin{equation}\label{decomposition}
M_{k,\nu}^{m}=N_{k,\nu}^{m}\oplus U_{k ,\nu}^{m},
\end{equation}
as $K[G]$-modules, where $N_{k,\nu}^{m}$ and $U_{k ,\nu}^{m}$ are generated respectively, by the sets
\begin{equation*}
\{ \omega_{i, \nu+e}^{m} | 0\leq i\lneqq k,\ 0\lneqq e\leq n(i,k,m)\} \textrm{ and } \{\theta_0 ,\ldots, \theta_k \},
\end{equation*}
 where
\begin{eqnarray}\label{theta's}
\theta_i &:=&\theta_{i}^{(k, \nu, m)}\nonumber\\
&=&\sum_{e=0}^{n(i,k, m)} b_{e}^{(i, k, m)}\omega_{i, \nu+e}^{m}=\begin{cases}  b_0 ^{(i, k, m)}\omega_{i,\nu}^m+\sum_{e=1}^{n(i,k, m)} b_{e}^{(i, k, m)}\omega_{i, \nu+e}^{m}, & \textrm{if }  i \neq k ,\cr
              \omega_{k,\nu}^m, & \textrm{if } i = k. \cr  \end{cases}
\end{eqnarray}
%Note that $N_{p^{n}-1,\nu}^{m}$ is generated by
%\begin{equation*}
%\{\omega_{k, \nu}^{m} | 0\leq k\lneqq p^n -1,\ 0\lneqq \nu \leq \Gamma_{k}(m)-2m\}.
%\end{equation*}
Notice that the decomposition in Eq. (\ref{decomposition}) is  a decomposition of 
$K$-vector spaces since we can select our model so that $b_0\neq 0$, see Remark \ref{rem24}.

The $K[G]$ action on $N_{k,\nu}^{m}$'s is given by
\begin{equation}\label{action on N}
\sigma_{\alpha}( \omega_{i, \nu+e}^{m})=\sum_{\delta=0}^{i}\sum_{\theta=0}^{n(\delta, i ,m)}\left(\begin{array}{clr}i\\ \delta \end{array}\right)\alpha^{i-\delta}b_{\theta}^{(\delta, i ,m)}\omega_{\delta, \nu+e+\theta}^{m}.
\end{equation}
Observing that 
\begin{eqnarray*}
0\lneqq e+\theta&\leq& n(i,k,m)+n(\delta, i, m)\\
&=&\Gamma_i (m)-\Gamma_{k}(m)+\Gamma_{\delta}(m)-\Gamma_i (m)\\
&=&n(\delta, k,m),
\end{eqnarray*}
we see that the action is well defined, i.e. $\omega_{\delta, \nu+e+\theta}^{m} \in N_{k,\nu}^{m}$. In a similar way, using Eq. (\ref{gaction on the basis}) for the $\theta_i$'s (compare also to \cite[p. 153, Eq. (7)]{csm}), we see  that  the $K[G$] action on $U_{k ,\nu}^{m}$ is given by
\begin{equation}\label{actionfinal}
\sigma_{\alpha}(\theta_i )=\sum_{l=0}^{i}\left(\begin{array}{clr}i\\ l \end{array}\right)\alpha^{i-l}\theta_l ,\textrm{ for }\  0\leq i \leq k.
\end{equation}
Thus the spaces $N_{k,\nu}^m$ and $U_{k ,\nu}^{m}$ are indeed $K[G]$ modules and Eq. (\ref{decomposition}),  is actually a $K[G]$-module decomposition of $M_{k,\nu}^{m}$.
\begin{remark} \label{rem24}
%\footnote{i added this remark in order to ensure that $ b_0 ^{(i, k, m)} \neq 0$}
In Eq. (\ref{theta's}), where we define the thetas, we can assume that $ b_0 ^{(i, k, m)} \neq 0$, in order to have the desired terms $\omega_{i,\nu}^m$. Indeed  $b_0 ^{(i, k, m)}$ is just the constant term of the polynomial $h_{i,k}(x)=\prod_{j=1}^r (x-a_j)^{\nu_{ji}(m)-\nu_{jk}(m)}$, which is zero when any of its roots $a_j$ equals to zero.  Thus we may assume, after a birational transformation  (i.e. an appropriate translation), that $a_j \neq 0$ in  Eq. (\ref{equation e.abelian}), for all $1\leq j\leq r$.
\end{remark}

We now present a counting argument:\\
Let $0\leq k_0 \leq p^n -1$ be maximal such that $\Gamma_{k_0}(  m) -2m \geq 0$ and set $0\leq \Gamma_{k_0}(  m) -2m= \nu_0 .$ Observe that with the above hypothesis $k_0 =p^n -1$. If we were in the case $m=1$ an admissible value for $k_0$ would be  $p^n -2$.
%\footnote{i added the $m=1$ case} 
\\
Recall that $ \Omega_L (m)$ is generated by
\begin{equation*}
%\bigcup_{i=0}^{k}
\{\omega_{k, \nu}^{m} | 0\leq k \leq k_0,\ 0\leq \nu \leq \Gamma_{k}(m)-2m\}.
\end{equation*}
\textit{Claim}:
\begin{equation}\label{dec0}
 \Omega_L (m)=N_{k_0,\nu_0}^{m}\bigoplus_{j=0}^{\nu_0} U_{k_0 ,j}^{m}, 
\end{equation}
as $K[G]$ modules, with 
$N_{k_0,\nu_0}$ generated by
\begin{equation}\label{N}
\{\omega_{k, \nu}^{m} | 0\leq k\lneqq k_0 ,\ 0\leq \nu_0\lneqq \nu \leq \Gamma_{k}(m)-2m\}
\end{equation}
and  for every $0\leq j\leq \nu_0$, each $ U_{k_0 ,j}^{m}$ is generated by
\begin{equation*}
\{ \theta_{i}^{(k_0 , j , m)}|\ 0\leq i\leq k_0 \} ,
\end{equation*}
where $ \theta_{i}^{(k_0 , j , m)}$'s are given by Eq. (\ref{theta's}) .  \\
\textit{Proof of the Claim}:\\
Note that $\dim_K U_{k_0 ,j}^{m}=k_0+1,$ for all $0\leq j \leq \nu_0.$

Remember, that $\theta_{i}^{(k_0, j, m)}$ equals to
\begin{eqnarray}
\theta_{i}^{(k_0, j, m)}&=& \sum_{e=0}^{n(i,k_0, m)} b_{e}^{(i, k_0, m)}\omega_{i,j+e}^{m}\nonumber \\
&=& 
\begin{cases}  b_0 ^{(i, k_0, m)}\omega_{i,j}^m+\sum_{e=1}^{n(i,k_0, m)} b_{e}^{(i, k_0, m)}\omega_{i, j+e}^{m}, & \textrm{if }  i \neq k_0 ,\cr
              \omega_{k_0,j}^m, & \textrm{if } i = k_0, \cr  \end{cases}\label{1}
%&&\textrm{  or alternatively }\nonumber \\
%&=&\begin{cases} \omega_{i,j+n(i,k_0, m)}^m+\sum_{e=0}^{n(i,k_0, m)-1} b_{e}^{(i, k_0, m)}\omega_{i, j+e}^{m}, & \textrm{if }  i \neq k_0 ,\cr
 %             \omega_{k_0,j}^m, & \textrm{if } i = k_0. \cr  \end{cases}\label{2}
\end{eqnarray}
We need the followings Propositions in order to prove the claim:
\begin{proposition}\label{claim1}
$U_{k_0, j}^m\cap N_{k_0 ,\nu_0 }^m=\{0\}, \textrm{ for all } 0\leq j\leq \nu_0.$
\end{proposition}
\begin{proof}
%Observe that the $K$-linear combinations of $\theta_{i}^{(k_0, j, m)}$'s for $0\leq j\leq \nu_0$ cannot cancel out, since every  $\theta_{i}^{(k_0, j, m)}$ is 
According to Eq. (\ref{1}),  every  $\theta_{i}^{(k_0, j, m)}$, would contain as a summand an $b_0 ^{(i, k_0, m)}\omega_{i,j}^{m}$, with $0\leq j \leq \nu_0$, but from the definition of $N_{k_0 ,\nu_0 }^m$,  Eq. (\ref{N}), these elements are not in $N_{k_0 ,\nu_0 }^m$. 
These elements $b_0 ^{(i, k_0, m)}\omega_{i,j}^{m}$ cannot be canceled out by linear 
combinations of elements in $U_{k_0,j}^m$. Therefore no linear combination of $\theta_{i}^{(k_0, j, m)}$ can be in $N_{k_0,\nu_0}^m$.
% 
% 
% Also every $\theta_{i}^{(k_0, j, m)}$ is a $K$-linear combination of distinct basis elements of $\Omega_L(m)$ which cannot cancel out. The same is true for the elements of $N_{k_0 ,\nu_0 }^m$. Thus  we cannot take  $\theta_{i}^{(k_0, j, m)}$ with $K$-linear  combinations of elements belonging to $N_{k_0 ,\nu_0 }^m$ as there is no way to get the summand $b_0 ^{(i, k_0, m)}\omega_{i,j}^{m}$, with $0\leq j \leq \nu_0$. On the other hand take an element $\omega_{k,\nu}^{m}\in  N_{k_0 ,\nu_0 }^m$ and observe that $\omega_{k,\nu}^{m}\neq \theta_{i}^{(k_0, j, m)}$, since $k\neq k_0$, nor can be written as a $K$-linear combination of $\theta_{i}^{(k_0, j, m)}$'s; every $\theta_{i}^{(k_0, j, m)}$ and every $K$-linear combination of $ \theta_{i}^{(k_0 , j , m)}$'s, with $ 0\leq i\leq k_0$ cannot cancel out in order to give $\omega_{k,\nu}^{m}$. Equivalently $\omega_{k,\nu}^{m}\notin U_{k_0 ,j}^{m}$ for all $j$'s .
\end{proof}
\begin{proposition}\label{claim2}
$U_{k_0, j}^m\cap U_{k_0, j'}^m=\{0\} $, for every $j\neq j'$, with $0\leq j\lneqq j'\leq \nu_0$.
\end{proposition}
\begin{proof}
Fix a $j$ and let $j'\neq j$. We may also assume that $j \lneqq j'$. 
We consider a linear combination of elements 
 $\theta_{i}^{(k_0, j', m)}$ in  $ U_{k_0, j'}^m$. 

Under the assumption $j \lneqq j'$, Eq. (\ref{1}) tells us that $b_0 ^{(i, k_0, m)}\omega_{i,j}^m$, a summand of a linear combination of $\theta_{i}^{(k_0, j, m)}$, is not a summand of a linear 
combination of the elements  $ \theta_{i}^{(k_0, j', m)}$, i.e $\omega_{i,j}^m\notin  \langle \theta_{i}^{(k_0, j', m)} \rangle_{0\leq i \leq k_0}$. 
% 
% On the other hand, using Eq. (\ref{2}) now, notice that every $\theta_{i}^{(k_0, j', m)}$ has a summand $ \omega_{i,j'+n(i,k_0, m)}^m$, but  $\omega_{i,j'+n(i,k_0, m)}^m \notin \langle \theta_{i}^{(k_0, j, m)}\rangle_{0\leq i\leq k_0}$, for all $j$.\\
\end{proof}

Finally, observe that the $M_{k_0, \nu_0}$ is a $\Omega_{L}(m)$ submodule of co-dimension $(k_0+1)\nu_0$ and, from Eq. (\ref{decomposition}) , the same is true for $N_{k_0,\nu_0}^{m}\oplus U_{k_0 ,\nu_0}^{m}$. The claim is then proved using Propositions \ref{claim1} and \ref{claim2}, if we notice that  $\dim_K\bigoplus_{j=0}^{\nu_0 -1 }U_{k_0, j}^m =\nu_0 (k_0 +1)$.

Observe also, using the Eq. (\ref {actionfinal}), that for every $j$ the $U_{k_0 , j}^{m} $'s that appear  in Eq. (\ref{dec0}), are $K[G]$ isomorphic. For example, using again Eq.  (\ref{theta's}), we can construct an isomorphism $f:U_{k_0 ,0}^{m}\longrightarrow U_{k_0 ,1}^{m}$ as follows; $f$  maps $\omega_{i, \nu}\mapsto \omega_{i,\nu+1}, $ with $0\leq \nu \leq \nu_0$ and $0\leq i\leq k_0$. So we can drop the $j$ subscript on the notation of  $ U_{k_0 ,j}^{m}$ and think of the $\bigoplus_{j=0}^{\nu_0} U_{k_0 ,j}^{m}$ as $\nu_0 +1 $ copies of  $ U_{k_0}^m$. Then, rewriting Eq. (\ref{dec0}) we have that
\begin{equation}\label{dec01}
 \Omega_L (m)=N_{k_0,\nu_0}^{m}\oplus [U_{k_0 }^{m}]^{\Gamma_{k_0}-2m+1}. 
\end{equation}

That finishes the zeroth step of the proof of the Theorem. Then we proceed to the first step. We take $N_{k_0,\nu_0}^{m}$, in place of $\Omega_F (m)$:
Let $0\leq k_1 \lneqq k_0= p^n -1$ be maximal such that $\Gamma_{k_1}(m) -2m \geq 0$ and set $0\leq \Gamma_{k_1}(m) -2m=\nu_1  .$ Then, repeating the claim in the previous step we, can see that
\begin{equation}\label{dec1}
N_{k_0,\nu_0}^{m}=N_{k_1,\nu_1}^{m}\bigoplus_{j=0}^{\nu_1} U_{k_1 ,j}^{m},
\end{equation}
where $
N_{k_1,\nu_1}^{m}$ is generated by
\begin{equation*}
\{ \omega_{k, \nu}^{m} | 0\leq k\lneqq k_1 ,\ 0\leq \nu_0 \lneqq \nu_1 \lneqq \nu \leq \Gamma_{k}(m) -2m\}
\end{equation*}
and  for every $0\leq j\leq \nu_1$, each $ U_{k_1 ,j}^{m}$ is generated by 
\begin{equation*}
\{\theta_{i}^{(k_1 , j , m)}|\ 0\leq i\leq k_1 \} ,
\end{equation*}
with $ \theta_{i}^{(k_1 , j , m)}$'s are given always by Eq. (\ref{theta's}).  Note that for $0\leq j\leq \nu_1$, $\dim_K U_{k_1 ,j}^{m}=k_1 +1$ and all the $U_{k_1 ,j}^{m}$'s  that appearing in Eq. (\ref{dec1}) are $K[G]$ isomorphic. There are exactly $\nu_1 -\nu_0 $ such modules, with $0\lneqq \nu_1-\nu_0 = \Gamma_{k_1}(m)-\Gamma_{k_0}(m)$. We can rewrite Eq. (\ref{dec1}), dropping the $j$ subscript, and thinking $\bigoplus_{j=0}^{\nu_1} U_{k_1 ,j}^{m}$ as $ \Gamma_{k_1}(m)-\Gamma_{k_0}(m)$ copies of $ U_{k_1}$:
\begin{equation}\label{dec11}
N_{k_0,\nu_0}^{m}=N_{k_1,\nu_1}^{m}\oplus [U_{k_1 }^{m}]^{ \Gamma_{k_1}(m)-\Gamma_{k_0}(m)}.
\end{equation}

Now we apply the above argument recursively to $N_{k_{\zeta},\nu_{\zeta}}^m$, for $0\leq \zeta \leq p^n -1$, continuing the above decomposition and replacing always the $N_{k_{\zeta},\nu_{\zeta}}^{m}$with the $N_{k_{\zeta-1},\nu_{\zeta-1}}^{m}$.

From Eq. (\ref{dec01}), (\ref{dec11}) and the repeated procedure,  we are now able to express $\Omega_F (m)$ as a direct sum of $U_{k_{\zeta },\nu_{\zeta}}^{m}$'s. Collecting these $K[G]$ modules of the same dimension we have 
%\begin{itemize}
%\item The $U_{k_{\zeta },\nu_{\zeta}}^{m}$'s that have dimension equal to $\dim_K U_{k_0 ,\nu_0}=p^n$,  have cardinal number equal to $\nu_0$, with $0\leq \nu_0\leq \Gamma_{k_0}(m)-2m$. So there are $ \Gamma_{p^n -1}(m)-2m+1$ of them.
%\item The $U_{k_{\zeta },\nu_{\zeta}}^{m}$'s that have dimension equal to  $\dim_KU_{k_1 ,\nu_1}=k_1 +1$, are exactly $\nu_1 -\nu_0 $, with $0\lneqq \nu_1-\nu_0\leq \Gamma_{k_1}(m)-\Gamma_{k_0}(m)$. Thus, there are $ \Gamma_{k_1}(m)-\Gamma_{k_0}(m)$ of them.
%\item $\cdots$
%\item The number of $U_{k_{\zeta },\nu_{\zeta}}^{m}$'s that have dimension over $K$ equal to $\dim_KU_{k_{p^n-1} ,\nu_{p^n-1}}=1$, is $ \Gamma_{0}(m)-\Gamma_{1}(m)$.
%\end{itemize}
%Equivalently
\begin{equation*}
\Omega_F (m)\simeq\oplus_{\zeta=0}^{p^n -1}[U_{k_{\zeta }}^{m}]^{\Lambda_{k_\zeta}},
\end{equation*}
with $\Lambda_{k_0}=\Gamma_{p^n -1}( m) -2m +1$ and $\Lambda_{k_\zeta}=\Gamma_{k_{\zeta}}(m)-\Gamma_{k_{\zeta-1}}(m)$, for all the steps: $1\leq \zeta \leq p^n -1$. From Eq. (\ref{actionfinal}) all $U_{k_{\zeta }}^{m}$ with the same dimension, say $j$, are $K[G]$-isomorphic, thus are isomorphic with $U_{j-1}^{m}$.  We re-index in order to be consistent with the dimension, letting $j-1=k_{\zeta}  $ (observe  that with this setting $k_{\zeta}+1=k_{\zeta-1}$), we obtain

\begin{eqnarray}\label{finaldec}
\Omega_F (m)&\simeq&\oplus_{j=1}^{p^n }[U_{j-1}^{m}]^{\Lambda_{j-1}}, \textrm{ or}\nonumber\\
\Omega_F (m)&\simeq &\oplus_{j=1}^{p^n }T_j .
\end{eqnarray}
The module $T_{p^n } $ is a direct sum of $\Gamma_{p^n -1}( m) -2m +1$ modules of dimension $p^n $ and $T_j $ is a direct sum of $\Gamma_{j-1}( m)-\Gamma_{j}(m)$ modules of dimension $j$, with $1\leq j \leq p^n -1$. 
%\begin{eqnarray*}
%&\bigcup_{i=0}^{k}&\{\omega_{i, \nu}^{m} | 0\leq i\leq k,\ 0\leq \nu \leq \Gamma_{k}(m)-2m\}\\
%=&\bigcup_{i=0}^{k_0}&\{\omega_{i, \nu_0}^{m} | 0\leq i\leq k_0 ,\ 0\leq \nu_0 \leq \Gamma_{k_0}(m)-2m\}\\
%=&\bigcup_{i=0}^{k_0-1}&\{\omega_{i, \nu_0}^{m} | 0\leq i\lneqq k_0 ,\ 0\leq \nu_0 \leq \Gamma_{k_0}(m)-2m\}\\
%&\oplus&\{ \theta_{i}^{(k_0 , \nu_0 , m)}|\ 0\leq i\leq k_0 \} .
%\end{eqnarray*}

We will now prove that the modules  $U_{j-1}^{m}$ are indeed indecomposable.\\
Let $\Omega_F (m) \simeq \oplus_{i=1}^{\eta}M_{i}$, be a decomposition in indecomposable $K[G]$-modules of the space of holomorphic $m$-(poly)differentials. Then 
\begin{equation}\label{lbound}
\eta \geq \sum_{j=1}^{p^n -1}\left(\Gamma_{j-1}( m)-\Gamma_{j}(m)\right)+\Gamma_{p^n-1}( m) -2m +1=\Gamma_0( m) -2m +1.
\end{equation}
Since $G$ is a $p$-group, we know that $\tau$, the one-dimensional trivial representation, is the only irreducible representation of $G$ (see \cite[p.187, Proposition 1.1]{Wein}). Then, if $M_i ^{G}$ denote the $K[G]$-submodule of fixed points of $M_i$, $M_i ^{G}$ would contain $\tau$ as a subrepresentation, so
%\footnote{What do you mean by $| \cdot |$? The module $M_i$ has infinitely many elements.}
 $\dim_K M_{i} ^{G}\geq 1$ . Thus
\begin{equation}\label{upbound}
 \eta \leq \dim_K \Omega_{F}^{G}(m).
\end{equation}
It is well known (see for example \cite[p. 271, V.2.2]{Farkas-Kra}, or \cite[p.83, Theorem III.4.6]{StiBo}), that differentials of $K(x)$ can be lifted to $G$-invariant differentials of $F$ via the Cotrace map. So $\Omega_{F}^{G}(m)=\{\kappa(dx)^{\otimes m} \ | \kappa \in K(x) \textrm{ with } \textrm{div}_F (\kappa(dx)^{\otimes m})\geq 0\}$ and we have
\begin{eqnarray}\label{Cotr}
\textrm{div}_F(\textrm{ Cotr}_{F/K(x)} (\kappa(dx)^{\otimes m}))&:=&\textrm{div}_F (\kappa(dx)^{\otimes m})\nonumber\\
&=&\textrm{Con}_{F/K(x)}(\textrm{div}_K (\kappa(dx)^{\otimes m})) +m\textrm{Diff}(F/K(x)).
\end{eqnarray}
Evaluating Eq. (\ref{Cotr}), for the places $P_i  \in \mathbb{P}_{F}$ we have for all $1\leq i \leq r$ and for all $\kappa(dx)^{\otimes m}\in \Omega_{F}^{G}(m)$, that 
\begin{equation*}
v_{P_i }(\textrm{div}_F (\kappa(dx)^{\otimes m}))= p^n v_{p_i} (\kappa) +m(\Phi(i) +1 )(p^n -1) \geq 0,
\end{equation*}
so 
\begin{eqnarray}\label{p-i}
 v_{p_i} (\kappa)& \geq& -\frac{m(\Phi(i) +1 )(p^n -1)}{p^n }, \textrm{or using Eq. (\ref{mi})}\nonumber \\
 v_{p_i} (\kappa)& \geq&-\nu_{i0} (m),\textrm{for all } 1 \leq i \leq r.
\end{eqnarray}
For $Q \in \mathbb{P}_{K(x)}$, with $Q\neq p_i , p_{\infty} $, taking $\kappa(dx)^{\otimes m}\in \Omega_{F}^{G}(m)$ and using Eq. (\ref{Cotr}) and (\ref{divdx^m}), we have that 
\begin{equation}\label{Q}
v_Q (\kappa)\geq 0,
\end{equation}
while for the infinite place of $K(x)$, the same hypothesis and Equations yield
\begin{equation}\label{p_infty}
v_{p_{\infty}}(\kappa )\geq 2m.
\end{equation}
Gathering Eq. (\ref{p-i}), (\ref{Q}) and (\ref{p_infty}), we can write $\Omega_{F}^{G}(m)$ in an alternative form, namely
\begin{equation*}
\Omega_{F}^{G}(m)=\left\{\frac{c(x)}{\prod_{i=1}^{r}(x-a_i )^{\nu_{i0}(m)}} |c(x) \in K[x], \deg c(x) \leq \sum_{i=1}^{r} \nu_{i0}(m)-2m\right\}.
\end{equation*}
Therefore $\dim_K \Omega_{F}^{G}(m)= \Gamma_{0}( m) -2m +1$. Using this fact, together with Eq. (\ref{lbound}) and (\ref{upbound}) we have that $\eta=\Gamma_{0}( m) -2m +1$. This shows that the $K[G]$-modules appearing to each $T_j $ in the decomposition of $\Omega_F (m)$ at Eq. (\ref{finaldec}), are all indecomposable for every $1\leq j\leq p^n $.

The theorem follows by 
letting  $W_j = U_{j-1}^{m}$ and $d_j = \Lambda_{j-1}$.
\end{proof}

%We know that rad$K[G]$ equals $I_G$, the augmentation ideal of $K[G]$. Moreover $K[G]$ is a local ring and its maximal left ideal is rad$K[G]$. Thus $K[G]/I_G $ is a simple and so an indecomposable $K[G]$ module. On the other hand the Trace (or some authors called it norm) submodule of $K[G]$  (i.e the submodule generated by the trace ) $N_G$, is also simple, thus an indecomposable $K[G]$ module. This fact is taken from \cite[p.3816, Proposition 6]{lara08}. More explicity,

\begin{remark}
%{\bf Na doume an einai telika swsth ayth h parathrhsh}
 %Observe from Eq. (\ref{actionfinal}) that 
%\begin{eqnarray*}
%\sigma_{\alpha}(\theta_i )&=&\sum_{l=0}^{i}\left(\begin{array}{clr}i\\ l \end{array}\right)\alpha^{i-l}\theta_l ,\textrm{ for }\  0\leq i \leq k\\
%\sigma_{\alpha}(\theta_i ) -\theta_i&=&\sum_{l=0}^{i-1}\left(\begin{array}{clr}i\\ l \end{array}\right)\alpha^{i-l}\theta_l 
%\end{eqnarray*}
%and when $i\neq 0$, the left hand is equal to $I_{G} U_{k,\nu}^m=\langle \sigma_{\alpha}-1 |\sigma_{\alpha}\in G\setminus \{\textrm{id.}\}, \alpha \in (\mathbb{F^*}_q , +)\rangle$, with $G\simeq(\mathbb{F}_q , +)$. This is the augmentation submodule of $U_{k,\nu}^m$. So 
%\begin{equation*}
 %I_{G} U_{i,\nu}^m= U_{i-1,\nu}^m, \textrm{ with } 0\lneq i\leq k
%\end{equation*}
%and notice that all the $W_j$, for $0\leq j\leq p^n$, arise in this way.\\
%On the other hand if $i=0$ then  $U_{0,\nu}^m$, which is generated by $\{\theta_0 \}$, is just the Trace $K[G]$ submodule of $U_{k,\nu}^m$, namely $N_{G}U_{k,\nu}^m=\langle \sum_{\alpha \in \mathbb{F}_q} \sigma_{\alpha} \rangle$, which is (the only) irreducible and hence an indecomposable $K[G]$ module (see  \cite[p.3816, Proposition 6]{lara08}).

Another way to see the indecomposability of the $W_j$'s, for $j\neq0$, is to notice that  if $W_j$, were decomposable, say $W_j=M_1\oplus M_2$, then each $M_i$ would contain a copy of $U_{0,\nu}^m$, and $W_j$ would contain $U_{0,\nu}^m \oplus U_{0,\nu}^m$ as a subrepresentation. But since $\dim_K W_j^G\leq1$, $(\dim_K W_{j}^G \leq \dim_K K[G]^G=1)$, that is  a contradiction.

\end{remark}

\section{A conjecture concerning abelian groups of order $p^n$.}
We strongly believe that for an arbitrary Galois $p$-extension, $F/E$, with abelian Galois group, one can calculate explicitly the $K[G]$-module structure of the space of holomorphic $m$-(poly)differentials:
\begin{conj}\label{conj}
Let $G$, be a $p$-group of automorphisms of $F$.
Set  $E=F^G$ and let  $g_E$ be the genus of $E$ and $g_F\geq 2$, the genus of $F$. Let  $m$ be  a natural number with $m>1$, $\delta_i$  the different exponent of the extension and $e_i$ the ramification indices of the ramified primes of $E$ in $F$. The regular representation of $G$ occurs $\Gamma_{p^{n}-1}(m)+(g_E -1)(2m-1)$  times in the representation of $G$ on $\Omega_F (m)$. For $k=1,\ldots,[F:E]-1$, the indecomposable representation of degree $k$ occurs $\Gamma_{k-1}(m) -\Gamma_{k}(m)$ times. Where $\Gamma_k (m)=\sum_{i } \nu_{ik}(m),$ are the Boseck invariants, with $i$ running over the ramified primes of $E$ in $F$ and the quantities  $ \nu_{ik}(m),$ are defined to be
$$
\left\lfloor\frac{m\delta_i + \{ \textrm{evaluation of the kth  E- basis element of F  by a normalized valuation of F }\} }{p^{e_i}}\right\rfloor,
$$
where the basis element is evaluated by a (normalized) valuation determined by a place of $F$ above a ramified place of $E$ and $\lfloor \cdot \rfloor$ denotes the integer part.
\end{conj}
For the case $m=1$ the above Conjecture has been proved in some cases, since has the form:
\begin{conj}\label{conj2} $\;$ \\
\begin{itemize}
 \item[(i)] 
Wild ramification:
If $G$ is a $p$-group, $m=1$ and there is a place in $\mathbb{P}_E$ that is totally ramified in $F$, with $F=E^G$, then the regular representation of $G$ occurs $d_{p^n}:=g_E$ times in the representation of $G$ on $\Omega_F (1)$. For $k=1,\ldots, p^n -1$, the indecomposable representation of degree $k$ occurs $d_{k}:=\Gamma_{k-1}(1) -\Gamma_{k}(1) +\beta$ times, where $\beta$ equals to $-1$, if $k=p^n-1$ and equals to  zero otherwise. The Boseck invariants are defined as before, by letting $m=1$.\\
\item[(ii)] 
Tame ramification:
If $m=1$ and $F/E$ is ramified of degree $n$, with $\textrm{g.c.d.}(n,p)=1$, then for $k=0,\ldots, n -1$, we have $n$ distinct  irreducible representations of degree $1$. The $k$th representation occurs $d_{n-k}:=\Gamma_{n-k}(1) -1+g_E+\beta$ times in the representation of $G$ in $\Omega_F(1)$, where $\beta$ equals to $1$, if $k=0$ and  zero otherwise.  The $\Gamma_k (1)=\sum_{i } \nu_{ik}(1),$ are the Boseck invariants, with $i$ running over the ramified primes of $E$ in $F$ and the quantities  $ \nu_{ik}(m),$ are defined to be
\begin{equation*}
\left\langle \frac{ \{ \textrm{evaluation of the kth  E- basis element of F  by a normalized valuation of F }\} }{e_i}\right\rangle,
\end{equation*}
where $e_i$ is the corresponding ramification index and $\left\langle \cdot \right\rangle $ denote the fractional part of $\cdot$.
\end{itemize}
\end{conj}
The first case is proved when $G$ is cyclic or an elementary abelian of order $p^n$. The second case is Hurwitz's Theorem  (see Theorem \ref{cyclic prime to p}). 
 Note that in all cases, Boseck invariants $\Gamma_{k}(m)$, defined to be the quantities that come out from Boseck's bases.

Let $\mathcal{C}$ denote the Cartier operator (see for example \cite[p. 349]{Salvador:00}), then from the theory of $\frac{1}{p}$--linear maps, it is well known that $\Omega_F (1)$ decomposes as
$$\Omega_F(1)=\Omega_{F}^s(1)\bigoplus\Omega_{F}^n (1),$$
where  $\Omega_{F}^s(1)$ denotes the semisimple part of $\Omega_F (1)$, that is the $K$ vector space spanned by the set 
$\{\omega \in \Omega_F (1) |\;  \mathcal{C}\omega =\omega\}$, and $\Omega_{F}^n (1)$  denotes the nilpotent part, the $K$ vector space spanned by 
$\{\omega \in \Omega_F (1) |\;  \mathcal{C}^i\omega =0, \textrm{ for some } i\geq1\}$. Now if the Conjecture \ref{conj2} (i) is proved then, coupled
 with the main result of  \cite{Salvador:00}, will give explicitly the structure of the nilpotent part of $\Omega_F (1)$, a problem that is open, as far as we know. Notice finally  
that we can calculate the   nilpotent part of $\Omega_F (1)$ in both the elementary abelian and the cyclic case, i.e. for the cases that this Conjecture 
has already been proved, combining the results of \cite{vm}, \cite{csm} and \cite{Salvador:00}.

The reason we believe that this conjecture is true is that we are able to prove it for the 
two extreme cases of abelian groups of order $p^n$, namely elementary abelian groups 
$\mathbb{Z}/p\mathbb{Z} \times \cdots \times \mathbb{Z}/p\mathbb{Z}$
and 
cyclic groups  $\mathbb{Z}/p^n \mathbb{Z}$. 
{ Since $G$ is abelian we can decompose it to a direct sum of 
% elementary abelian $p$ groups and 
cyclic $p$-groups. Thus  a function field 
having $G$ as its Galois group can result as the compositum of (some) cyclic function fields found in subsection \ref{Galoismodulecyclic}
with  elementary abelian function fields from section \ref{elab}. A difficult that now arises for such extensions is that we do not know
how Boseck invariants, and hence bases (even in the $m=1$ case) behave under such compositums.}
\label{consection}
\section{An application to local deformation functors}
\label{4section}
Let $G$ be a $p$-group.
It was observed in \cite{KoJPAA06} that the tangent space of the global deformation functor 
$H^1(G,\mathcal{T}_X)$
can be computed in terms of covariants of $2$-holomorphic differentials by 
\begin{equation} \label{GDF}
 H^1(G,\mathcal{T}_X)=\Omega_X^{\otimes 2}\otimes_{K[G]}K.
\end{equation}
In this section we will use the results we obtained so far in order to express the dimension 
of the above spaces in terms of the Boseck invariants. We will use the global deformation 
functor approach in order to study the tangent space $H^1(G,\mathcal{T}_{K[[t]]})$ of the  local deformation functor in the 
sense of J.Bertin and A. M\'ezard \cite{Be-Me}. 
This can be done by considering Katz-Gabber covers \cite{KaGa} of the projective line, i.e.
Galois cover $\pi:X \rightarrow \mathbb{P}^1$ with only one full ramification point and Galois group $G$.
For first order infinitesimal deformations of the curve $X$ with the automorphism group $G$, there 
is a splitting of the tangent space $H^1(G,\mathcal{T}_X)$:
\begin{equation} \label{eee1}
 H^1(G,\mathcal{T}_X)=H^1(X/G,\pi_*^G(\mathcal{T}_X))\oplus H^1(G,\mathcal{T}_{K[[t]]}).
\end{equation}
For the dimension of the space $H^1(X/G,\pi_*^G(\mathcal{T}_X))$ we have an explicit formula, namely
\begin{equation} \label{eee2}
 \dim_K H^1(X/G,\pi_*^G(\mathcal{T}_X))=3 g_{X/G}-3 + \lc \frac{\delta}{p^n} \rc,
\end{equation}
where $\delta$ is the local contribution to the different at the unique ramification point
\cite[Eq. (38)]{KontoANT}.

{\bf Case 1.} The group $G$ is cyclic.
In this case each of  the indecomposable components of Theorem \ref{basic thrm} has an one dimensional 
covariant subspace thus 
\begin{eqnarray} \label{eee3}
 \dim_K \Omega_X^{\otimes 2}\otimes_{K[G]}K &= &\sum_{\nu=1}^{p^n}  d_{\nu} \\
 & =&3(g_{X/G}-1)+\Gamma_0(2)=-3+\lf \frac{2 \delta}{p^n} \rf. \nonumber
\end{eqnarray}
Using Eq. (\ref{eee1}),(\ref{eee2}), (\ref{eee3}) we arrive at 
\[
 \dim_K H^1(G,\mathcal{T}_{K[[t]]})=\lf \frac{2 \delta}{p^n} \rf -  \lc \frac{\delta}{p^n} \rc,
\]
which coincides with the computation of \cite[Prop. 4.1.1]{Be-Me}.

{\bf Case 2.} The group $G$ is elementary abelian.
In this case we will use Theorem \ref{fff} in order to arrive to 
\begin{equation} \label{mmmm}
 \dim_K H^1(G,\mathcal{T}_{K[[t]]})=\sum_{j=1}^{p^n} d_j\cdot \dim_K (W_j \otimes _{K[G]}K )+3 -
\lc \frac{\delta}{p^n} \rc.
\end{equation}
\begin{proposition}
 For the dimension $W_j \otimes _{K[G]}K$ we compute 
\[
 \dim_K W_j \otimes _{K[G]}K=\left\{
\begin{array}{ll}
 1 & \mbox{ if }  1 \leq j \leq p \\
 2 & \mbox{ if }  p+1\leq j \leq p^n  \mbox{ and } j \neq 0 \mbox{\textup{mod}} p \\
 1 & \mbox{ if } j=0 \mbox{\textup{mod}} p
\end{array}
\right.
\]
\end{proposition}
\begin{proof}
We identify an elementary abelian group of order $p^n$ with the additive group of the 
field $\mathbb{F}_{p^n}$. The field $\mathbb{F}_{p^n}$ is an $\mathbb{F}_p$ vector space 
with basis $1,e,e^2,\ldots,e^{n-1}$ for some element $e$.
Every element $a\in \mathbb{F}_{p^n}$ gives rise to an automorphism $\sigma_a$.
We will denote by $\bar{W}_j=W_j \otimes_{K[G]}K$.
The modules $\bar{W}_j$ are given by $W_j/(\sigma_a (w)-w)$, where $a$ runs over $\mathbb{F}_{p^n}$ and 
$w$ runs over $W_j$.

The module $W_1=\langle \theta_0\rangle$ and is already $G$-invariant.
Observe  that $W_2=\langle \theta_0,\theta_1 \rangle$ and the action is given 
by $\sigma_a(\theta_0)=\theta_0$, $\sigma_a(\theta_1)=\theta_1+a \theta_0$.
Therefore we have only one relation in the module of covariants $\bar{W}_2$
namely $\sigma_a(\theta_1)-\theta_1=a \theta_0$ which implies that the image 
$\bar{\theta}_0$ in $\bar{W}_2$ is zero. 

 The module $W_3$ is generated by $\theta_0,\theta_1,\theta_2$ and the 
relation $\bar{\theta}_0=0$ is inherited in $\bar{W}_3$. We also have the relation 
$\sigma_a(\theta_2)=\theta_2+ a \theta_1 + a^2 \theta_0$, which implies that $\bar{\theta}_1=0$ 
in $\bar{W}_3$.

{We proceed by induction. For the inductive step we assume that for $j+1 \leq p$ we have 
the relation $\bar{\theta}_0,\ldots,\bar{\theta}_{j-2}=0$ in $\bar{W}_j$. Then 
$\sigma_a(\theta_{j})=\theta_{j} + a \theta_{j-1} + L$ where $L$ is an $\mathbb{F}_{p^n}$ linear combination of 
$\theta_\nu$ with $\nu \leq j-2$ that have zero image in the module of covariants.
Thus $\bar{\theta}_{j-1}=0$ in $\bar{W}_{j+1}$ and $\bar{W}_{j+1}=\langle \bar{\theta}_j \rangle$.

For the module $W_{p+1}$ the situation changes: We have $\sigma_a (\theta_p)=\theta_p+ a^p \theta_0$, which 
does not give any new relation. Therefore $\bar{W}_{p+1}$ is two dimensional generated by 
$\bar{\theta}_{p-1},\bar{\theta}_p$. 

We proceed by induction. For the modules $\bar{W}_{p+\nu+1}$, $1 \leq  \nu < p-1$ the inductive hypothesis is that 
$\bar{W}_{p+\nu+1}=\langle \bar{\theta}_{p-1},\bar{\theta}_{p+\nu}\rangle$.
We compute that for $1<\nu \leq  p-1$
\begin{equation}\label{eisenbud}
 \sigma_a(\theta_{p+\nu})=\sum_{\mu=1}^\nu \binom{p+\nu}{p+\mu} a^{\nu-\mu}\theta_{p+\mu}
 + \binom{p+\nu}{p-1}a^{\nu+1}\theta_{p-1}. 
\end{equation}
Notice that according to \cite[prop. 15.21]{Eisenbud:95} $\binom{p+\nu}{p-1}=0$ unless $\nu=p-1$.
Therefore for $\nu < p-2$ and after some computations, we arrive at
\[
 \sigma_a(\theta_{p+\nu+1})=\theta_{p+\nu+1}+ \nu a \theta_{p+\nu},
\]
and this implies that $\bar{\theta}_{p+\nu}=0$ in $\bar{W}_{p+\nu+2}$, thus 
$\bar{W}_{p+\nu+2}=\langle \bar{\theta}_{p-1},\bar{\theta}_{p+\nu+1}\rangle$.

Since $\bar{W}_{2p-1}=\langle \bar{\theta}_{p-1}, \bar{\theta}_{2p-2} \rangle$, we can now compute from Eq. (\ref{eisenbud}) for $\nu=p-1$, 
\[\sigma(\theta_{2p-1})=\theta_{2p-1}-a \theta_{2p-2} + a^p \theta_{p-1}.\]
Therefore in $\bar{W}_{2p}$ we have the relations 
$a \bar{\theta}_{2p-2} = a^{p} \bar{\theta}_{p-1}$ for $a \in \mathbb{F}_{p^n}$.
Taking $a=1$ we obtain $\bar{\theta}_{2p-2}=\bar{\theta}_{p-1}$ and then by taking $a=e$ we 
have $ \bar{\theta}_{2p-2} =e^{p-1} \bar{\theta}_{2p-2}$, therefore $\bar{\theta}_{p-1}=\bar{\theta}_{2p-2}=0$.
Thus $\bar{W}_{2p}=\langle \bar{\theta}_{2p-1} \rangle$.

We now continue to $\bar{W}_{2p+1}$ by computing  that 
\[
 \sigma_a (\theta_{2p})=\theta_{2p} +2 a^p \theta_p + a^{p^2} \theta_0,
\]
thus $\bar{W}_{2p+1}=\langle \bar{\theta}_{2p-1},\bar{\theta}_{2p} \rangle$ is of dimension $2$.

If $n \geq 2$ we proceed the same way: The modules $\bar{W}_{2p+\nu}$, $1\leq \nu \leq p-1$ are 2-dimensional
and $\bar{W}_{3p}$ is one dimensional. A final inductive argument shows that the $\bar{W}_{\lambda p +\nu}$, for $\lambda\leq p^{n-1}$ have 
the desired dimensions.}

\end{proof}
We now can give a closed formula for the sum given in eq. (\ref{mmmm}).
By the construction of the Katz-Gabber cover there is only one ramified point and we are interested for 2-holomorphic 
differentials ($m=2$) so we set 
\[
\Gamma_j (2)=\Gamma_j :=\lf \frac{2(p^n-1)(\Phi+1)+j \Phi}{p^n} \rf. 
\]
We compute:
\[
\sum_{j=1}^{p^n} d_j\cdot \dim_K (W_j \otimes _{K[G]}K )
=\]
\[
=\sum_{j=1}^{p-1} d_j\cdot \dim_K (W_j \otimes _{K[G]}K )+ 
\sum_{j=p}^{p^n} d_j\cdot \dim_K (W_j \otimes _{K[G]}K )=
\]
\[
 =\Gamma_0-\Gamma_{p-1} + 2(\Gamma_{p-1} -3) -\sum_{\nu=0 \atop p\mid \nu}^{p^n-1} \left( \Gamma_{\nu-1} -\Gamma_{\nu} \right)
=\]
\begin{equation} \label{ffff}
 =\Gamma_0 +\Gamma_{p-1} -6 -\sum_{\nu=0 \atop p\mid \nu}^{p^n-1} \left( \Gamma_{\nu-1} -\Gamma_{\nu} \right).
\end{equation}
This formula should give the same results with  the formula given in
\cite{KontoANT}. 
Giving a direct proof that the two formulas coincide is a complicated task to do. 
However using the Magma algebra system we checked that Eq. (\ref{ffff}) 
coincides with the formula given in \cite{KontoANT}, for all   choices of $G,p$ that we tried.
The magma \cite{magma}  program used to 
compute them is available at \texttt{http://myria.math.aegean.gr/\~{ }kontogar/sk/}.
% \begin{center} 
% % use packages: array
% \begin{tabular}{|l|l|l|} 
% \hline
% Cases  & Dimensions &  $\dim_K H^1(G,\mathcal{T}_X)$ \\
% \hline
% $p=7,n=2$ & 
% $\begin{array}{l|llllll}
% j&  8 & 18 & 28 & 37 & 47 & 49 \\
% \hline
% d_j &1 & 1 & 1 & 1 & 1 & 3  \\
% \hline
% \delta_j &    2 & 2 & 1 & 2 & 2 & 1 
%              \end{array}
% $
% & $12$
%  \\ 
% \hline
% $p=7,n=3$ & 
% $\begin{array}{l|llllll}
% j&  67 & 135 & 204 & 273 & 341 & 343 \\
% \hline
% d_j &1 & 1 & 1 & 1 & 1 & 3  \\
% \hline
% \delta_j &    3 & 3 & 3 & 2 & 3 & 1 
%              \end{array}
% $
%  & 
% $17 $
% \\ 
% \hline
% $p=13,n=2$ & 
% $\begin{array}{l|llllll}
% j&  32 & 66 & 100& 133 & 167 & 169  \\
% \hline
% d_j &1 & 1 &1  & 1 & 1  & 3  \\
% \hline
% \delta_j &    2 &2& 2 & 2 & 2  & 1 
%              \end{array}
% $
% & $13$\\
% \hline
% \end{tabular}
% \end{center}

{\bf Remark:} If the conjectures given in section \ref{consection} are proved then 
we have a method to compute the tangent space for the deformation space of curves 
with automorphism in the case of abelian groups.   

\textbf{Acknowledgment.} The author would like to thank %Prof. 
A. Kontogeorgis for all his valuable comments  while reading earlier versions of this paper, for the enlightening conversations we had, and for being a wonderful companion during this first journey to mathematical research. In section \ref{4section},  the work needed in order to confirm our results is completely due to him.

\def\cprime{$'$}
\providecommand{\bysame}{\leavevmode\hbox to3em{\hrulefill}\thinspace}
\providecommand{\MR}{\relax\ifhmode\unskip\space\fi MR }
\providecommand{\MRhref}[2]{%
  \href{http://www.ams.org/mathscinet-getitem?mr=#1}{#2}
}
\providecommand{\href}[2]{#2}

\def\cprime{$'$}
\providecommand{\bysame}{\leavevmode\hbox to3em{\hrulefill}\thinspace}
\providecommand{\MR}{\relax\ifhmode\unskip\space\fi MR }
\providecommand{\MRhref}[2]{%
  \href{http://www.ams.org/mathscinet-getitem?mr=#1}{#2}
}
\providecommand{\href}[2]{#2}

\end{document}